\documentclass[preprint]{imsart}

\RequirePackage{amsthm,amsmath,amsfonts,amssymb}
\RequirePackage[numbers]{natbib}
\RequirePackage{mathtools,dsfont}
\RequirePackage[colorlinks,citecolor=blue,urlcolor=blue]{hyperref}
\RequirePackage{graphicx}
\RequirePackage[normalem]{ulem}
\RequirePackage[dvipsnames]{xcolor}

\startlocaldefs
\numberwithin{equation}{section}
\theoremstyle{plain}
\newtheorem{theorem}{Theorem}[section]
\newtheorem{example}[theorem]{Example}
\newtheorem{lemma}[theorem]{Lemma}
\newtheorem{proposition}[theorem]{Proposition}
\newtheorem{cor}[theorem]{Corollary}
\theoremstyle{remark}
\newtheorem{remark}[theorem]{Remark}
\newtheorem{assumption}[theorem]{Assumption}
\DeclarePairedDelimiter{\floor}{\lfloor}{\rfloor}

\newcommand{\N}{\mathbb{N}}
\newcommand{\R}{\mathbb{R}}
\newcommand{\Q}{\mathbb{Q}}

\newcommand{\Z}{\mathbb{Z}}

\newcommand{\p}{\mathbb{P}}

\newcommand{\cA}{\mathcal{A}}

\newcommand{\cC}{\mathcal{C}}

\newcommand{\cE}{\mathcal{E}}

\newcommand{\cF}{\mathcal{F}}

\newcommand{\cL}{\mathcal{L}}

\newcommand{\cX}{\mathcal{X}}

\newcommand{\E}{\mathbb{E}}
\newcommand{\Cov}{\mathrm{Cov}}
 
\renewcommand{\phi}{\varphi}
\renewcommand{\epsilon}{\varepsilon}
\newcommand{\diff}{\mathrm{d}}

\newcommand{\diam}{\operatorname{diam}}

\newcommand{\ol}{\overline}
\newcommand{\wt}{\widetilde}
\newcommand{\wh}{\widehat}
\newcommand{\1}[1]{\,\mathds{1}\! \left\{ #1 \right\} }

\newcommand{\Filt}{\mathbb{F}\mathrm{ilt}}

\newcommand{\cmX}{\cX^{(m)}}
\newcommand{\pd}{\operatorname{PD}}

\endlocaldefs

\allowdisplaybreaks

\begin{document}

\begin{frontmatter}
\title{
On the stability of the filtration functions for dependent data with applications to break detection
}
\runtitle{Stability of filtration functions and break detection}
\begin{aug}
\author{\fnms{Johannes} \snm{Krebs}  
\ead[label=e1]{johannes.krebs@ku.de}}
\address{KU Eichstätt, Ostenstraße 28, 85072 Eichstätt, Germany.
 }

\author{\fnms{Daniel} \snm{Rademacher}  
\ead[label=e2]{d.rademacher@tugraz.at}}
\address{TU Graz, Kopernikusgasse 24/III, 8010 Graz, Austria.
}

\runauthor{Krebs and Rademacher}
\affiliation{KU Eichstätt and TU Graz}

\end{aug}

\begin{abstract}
In this paper, we study the stability of commonly used filtration functions in topological data analysis under small perturbations of the underlying nonrandom point cloud.
Relying on these stability results, we then develop a test procedure to detect and determine structural breaks in a sequence of topological data objects obtained from weakly dependent data. 
The proposed method applies for instance to statistics of persistence diagrams of $\R^d$-valued Bernoulli shift systems under the \v Cech or Vietoris-Rips filtration.
\end{abstract}

%

%
%
\end{frontmatter}

\section{Introduction}
We study stability properties of filtration functions in topological data analysis (TDA) and building on these, we construct a framework for  tests of break detection in a series of statistics obtained from topological objects.
This article requires an intermediate knowledge of TDA, nevertheless, we introduce the main objects. We refer to the introduction by Chazal and Bertrand \cite{chazal2021introduction} and the monograph of Boissonnat et al. \cite{boissonnat2018geometric} for a nice introduction to topological methods in the context of data analysis.

The manuscript consists of three parts. The first part is ``nonrandom'' and considers stability results of filtration functions of the following type: Let $\phi$ be a map which returns the filtration time $\phi[J](x)$ of the simplex $x_J=(x_j:j\in J)^T$ for a finite set $J \subseteq [r] \coloneqq \{1,\ldots,r\}$ and a vector $x=(x_1,\ldots,x_r)^T$. Each $x_i$ is an element of $M \subset\R^d$ which is a compact, convex and full $d$-dimensional set in our setting. Then the entire filtration function $\phi$ takes $L$ different values 
$$
	\phi[J_1](x) < \phi[J_2](x) < \ldots < \phi[J_L](x)
$$
for certain minimal subsets $J_1 \subset J_2 \subset \ldots \subset J_L \subseteq [r]$ (which are characterized below in more detail). We also write $L^x$ and $J_i^x$ (for $1\le i\le L^x$) to highlight the dependence on $x$ of these minimal sets. We study the functional
\begin{align}\begin{split}\label{D:RhoFunction}
	\rho\colon M^r\to [0,\infty), \quad &x\mapsto \sup\Big\{ \epsilon \ge 0 \ | \ \text{for all } y\in B^{\otimes r}(x;\epsilon) :\\
	&\qquad\qquad\quad L^x = L^y \text{ and } J_i^x = J_i^y \text{ for } 1\le i \le L^x \Big\},
	\end{split}
\end{align}
where $B^{\otimes r}(x;\epsilon) = \prod_{i=1}^r B(x_i; \epsilon)$ is the product of the closed $\epsilon$-neighborhoods w.r.t.\ the Euclidean 2-norm. Here we build on contributions of Chazal and Divol \cite{chazal2018density}.
In our setting the functional $\rho$ is Lebesgue-a.e. positive on $M^r$ and allows us to quantify the stability of the filtration under sufficiently small perturbations of the point set as follows: We give sufficient conditions for the existence of an $\alpha \in (0,1]$ for which the $(d\cdot r)$-dimensional Lebesgue-measure satisfies 
\begin{align}\label{E:Rho-Assumption}
	\lambda^{d r}( \{ x\in M^r \mid \rho(x) \leq t\} ) \le C^* t^{\alpha}, \qquad t\in[0,1].
\end{align}
The constant $C^*$ only depends on $r,d$ and $M$. We show that the \v Cech and Vietoris-Rips filtration function attain the exponent $\alpha=1$ and satisfy
$$
			\lambda^{d r} \big( \{x\in M^r: \rho(x) \le t \} \big) \le 
			\begin{cases}
				 c r^4 t, & \text{VR-filtration,}\\
				c r^{2(d+1)} t, & \text{\v Cech filtration}
			\end{cases}
$$ 
for a constant $c$ which only depends on $M$ and $d$ but neither on $r$ nor $t$.

This stability result is then applied in the second part to quantify the dependence of persistence diagrams $(\pd_{k,t})_t$ of dimension $k$ which are obtained from weakly dependent point clouds. In more detail, we study a sequence of topological statistics
\begin{align}\label{E:Intro1}	
	f(Z_{k,1}), \ldots, f(Z_{k,n})
\end{align}
which are derived from the feature vector representation $Z_{k,1},\ldots,Z_{k,n}$ of the sequence $(\pd_{k,t})_t$ by means of a Lipschitz continuous function $f$. Here the $Z_{k,t}$ encode one-to-one the $k$-dimensional persistence diagram $\pd_{k,t}$ of a point cloud $\cX_t = (X_{t,1},\ldots,X_{t,r}) \subset M^r$.
We show that the $(Z_{k,t})_t$ are $L^p$-$m$-approximable in the spirit of Hörmann and Kokoszka \cite{hormann2010weakly} provided the $(\cX_t)_t$ are $L^{p'}$-$m$-approximable for some higher order $p'$ which also involves the coefficient $\alpha$ from above.
Finally in the third part, we develop a CUSUM statistic-based test procedure for the data in \eqref{E:Intro1} in the spirit of Aue et al. \cite{aue2009break} for the statistics
\begin{align*}
	H_0\colon \pd_{k,1} = \ldots = \pd_{k,n},
\end{align*}
against the alternative 
\begin{align*}
	H_1\colon \pd_{k,1} = \ldots = \pd_{k,v^*} \not= \pd_{k,v^*+1} = \ldots = \pd_{k,n}.
\end{align*}
Here $v^*=\floor{n \theta}$ is a change point determined by an unknown $\theta\in (0,1)$. Under the null hypothesis, the limit distributions are known (see for instance \cite{aue2009break}). Under the alternative and in the specific case of a change in mean of the statistics $(f(Z_{k,t}))_t$, we show that the CUSUM statistic-based change point estimator $\wh \theta_n$ is strongly consistent with a rate $n^{-1} a_n$, where under reasonable conditions $(a_n)_n$ only has to satisfy $a_n / (n^{1/2} \log\log n) \to\infty$.\\

The theory on topological methods on time series or dependent data is rather sparse. Seversky et al.\ \cite{seversky2016time} and Umeda \cite{umeda2017time} were among the first who propose topological methods in the context of time series. Recent contributions related to topolical inference for time series data focus on dependence given by $\beta$-mixing conditions, see for instance Reise et al. \cite{reise2023topological} and Chazal et al. \cite{chazal2024topological}.
As already pointed out, our methods are developed under the assumption of $L^p$-$m$-approximable time series, a concept which was popularized by contributions around 2010 (see e.g. H{\"o}rmann and Kokoszka \cite{hormann2010weakly} and Aue et al. \cite{aue2009break}); related conditions however date back at least to Ibragimov \cite{ibragimov1962some}. To be more specific, we assume that the underlying point cloud process $\cX_t$ has a Bernoulli shift structure of the form $F(\epsilon_t,\epsilon_{t-1},\ldots)$, where $(\epsilon_t)_t$ are certain iid elements aggregated suitably by a function $F$. This approach of physical dependence has successfully been employed by Wu \cite{wu2005bahadur, wu2005nonlinear} or Berkes et al. \cite{berkes2014komlos}. The $L^p$-$m$-framework contains broad classes of standard time series models such as (multivariate) ARMA and GARCH processes or certain Markov chains (see also Rosenblatt \cite{rosenblatt2009stationary} and the references therein) and a large variety of limit theorems exists for this type of dependent data. A related concept of weak dependence which relies solely on projective methods has been formulated by Dedecker et al.\ \cite{dedecker2007weak}.

Islambekov et al.\ \cite{islambekov2020harnessing} and Zheng et al. \cite{zheng2023percept} examined change point detection based on topological methods using data experiments.
However, a consistent theoretical framework for break detection based on topological methods is still missing to the best of our knowledge.
The detection of structural changes in a univariate or multivariate series of observations is one of the fundamental problems in mathematical statistics and has enjoyed quite some popularity ever since. 
For statistical tests on structural breaks in various settings, we refer to the contributions of James et al. \cite{james1987tests}, Carlstein et al. \cite{carlstein1994change}, Aue et al. \cite{aue2009break,aue2013structural,aue2018detecting}, Matteson and James \cite{matteson2014nonparametric}, Chen and Zhang \cite{chen2015graph}, Jirak \cite{jirak2015uniform}, Chen \cite{chen2019sequential}, Wang and Samworth \cite{wang2018high}, Wang et al. \cite{wang2021optimal}.
Recently, tests for so-called relevant differences have appeared in the literature, see the contributions of
Dette and Wied \cite{dette2016detecting}, Dette and Wu \cite{dette2019detecting} and Aue et al.\ \cite{aue2019two}.
Moreover, owed to increased complexity of data, detection of change points in general Fr{\' e}chet objects have gained popularity, too; see for instance Dubey and Müller \cite{dubey2019frechet, dubey2020frechet}.\\

We shortly give the main abbreviations used throughout the paper. $\|\cdot\|_q$ is the $q$-norm on the Euclidean space. In particular, we set $\|\cdot\| \coloneqq \|\cdot\|_2$ for the 2-norm and $B(x;\epsilon)$ is the closed $\epsilon $-neighborhood of $x$ w.r.t.\ $\|\cdot\|$.
Given an index set $I$ and a sequence $(x_i: i\in I)$ of real numbers or random variables and $A\subset I$, we set $x_A = (x_i: i\in A)$. If $A\subseteq\R^d$ and $t>0$, we set $A^{(t)} = \{x\in\R^d: d(x,A)\le t\}$ and $A^{(-t)} = \{x\in A: d(x,\R^d\setminus A) \ge t \}$. We denote the smallest closed bounding ball (not bounding sphere) of a generic $x_{J}$, $J \subseteq [r]$, 
by $B^*(x_{J}) \coloneqq B(z(x_J), R(x_J))$ with circumcenter  $z(x_{J})\in \R^d$ and circumradius $R(x_{J})>0$, i.e. $x_{J} \subset B^*(x_{J})$ and $x_{ J} \not\subset B(z',R')$ for all $R'< R(x_{ J})$ and all $z'\in\R^d$.
Let $(X_n)_n$ and $(Z_n)_n$ be sequences of real-valued random variables. We write $X_n = o_{a.s.}(|Z_n|)$ if $X_n/|Z_n|\to 0$ a.s. We denote the weak convergence of random variables by the symbol $\Rightarrow$.\\

We give a short outline of the remainder of the paper. In Section~\ref{Section_Stability} we give the main results for the stability of popular filtration functions. These findings are then applied in Section~\ref{Section_FeatureVectors}. Here we study the dependence and asymptotic behavior of a sequence of statistics of persistence diagrams obtained from weakly dependent point clouds. Building on the results from that section, we then develop a framework to assess breaks in a sequence of weakly dependent persistence diagrams in Section~\ref{Section_Test}. Finally, we give the mathematical details of all statements in Section~\ref{Section_Proofs}.

\section{Stability of filtrations under small perturbations}\label{Section_Stability}

The results in this section apply to a larger class of filtration functions and not only to the \v Cech and Vietoris-Rips filtration. To this end we introduce some general notation.

Let $M$ be a compact, convex and full $d$-dimensional subset of $\R^d$.
Given a point set $x=(x_1,\ldots,x_r) \in M^r$ we define a filtration $\Filt(x) = \{ \Filt(x;s)\mid s\geq 0 \}$ by using a
filtration function $\phi$ as follows:
Let $\mathcal{F}_r$ denote the collection of all $2^r-1$ nonempty subsets of $[r]$. Suppose 
$\phi = (\phi[J])_{J\in\mathcal{F}_r} \colon M^r \to \R^{2^r-1}$ is a continuous function. Then a $k$-simplex $x_J=\{ x_{i_0},\ldots, x_{i_k} \}$, $J=(i_0,\ldots,i_k)\in \mathcal{F}_r$ is added to the filtration at time $\phi[J](x)$, i.e., 
\begin{align*}
	x(J) \in \Filt(x;s) \quad:\Longleftrightarrow\quad \phi[J](x) \leq s \qquad \forall J\in \mathcal{F}_r.
\end{align*} 

We impose some regularity conditions on $\phi$ and require the conditions (K1)-(K5), resp., conditions (K1)-(K4) and (K5') of Chazal and Divol \cite{divol2019density}. We cite these conditions below (in the original framework) and refer to Chazal and Divol for a more detailed discussion and some background on subanalytic geometry.
\begin{itemize}
	\item [(K1)] For each $J\in\mathcal{F}_r$, $\phi[J](x)$ only depends on $x(J)$.
	\item [(K2)] For each $J\in \mathcal{F}_r$ and permutation $\pi$ of $J$ we have for $x\in M^r$ and $x_\pi = (x_{\pi(1)},\ldots,x_{\pi(r)})$ that $\phi[J](x) = \phi[J](x_\pi)$.
	\item [(K3)] If $J\subseteq J'$, then $\phi[J] \le \phi[J']$.
	\item [(K4)] For $J\in \mathcal{F}_r$ and $j\in J$, if $\phi[J](x)$ is not a function of $x_j$ on some open set $U$ of $M^r$, then $\phi[J] = \phi[J\setminus \{j\}]$ on $U$.
	\item [(K5)] The function $\phi$ is subanalytic and the (partial) gradient of each of its entries is nonvanishing a.s.e.
	\item [(K5')] The function $\phi$ is subanalytic and the (partial) gradient of its entries $J$ of size greater than 1 is nonvanishing a.s.e. Moreover, $\phi[\{j\}] \equiv 0$ for each singleton $\{j\}$.
\end{itemize}
In the following, we will mostly discuss the \v Cech and Vietoris-Rips filtration as standard choices.
The \v Cech filtration function is
$$
	\phi[J](x) = \inf_{\epsilon >0} \Big\{	\bigcap_{j\in J} B(x_j,\epsilon)\neq \emptyset \Big\} \qquad (J\in\cF_r, \quad x\in M^r).
$$
Hence, the filtration time $\phi[J](x)$ is the radius of the circumsphere of the set $\{x_j: j\in J\}$, viz. the radius of the smallest $d$-dimensional sphere that contains all $x_j$. Note that depending on the size of $J$ and the location of the $x_j$ not all elements $x_j$ lie necessarily on the circumsphere.

The Vietoris-Rips filtration function is
$$
	\phi[J](x) = \max_{i,j\in J} \|x_i-x_j\| \qquad (J\in\cF_r, \quad x\in M^r).
$$
It is well known that the \v Cech and Vietoris-Rips filtration both satisfy (K1) to (K4) and (K5'), see Chazal and Divol \cite{divol2019density}. 
In the following, we define 
$$
	T = \sup\{ \phi[J](x) : J\subseteq [r], x\in M^r \} < \infty.
	$$
Then the \v Cech filtration satisfies $2T \le \diam M$ and the Vietoris-Rips filtration $T =  \diam M$.

For a given $x \in M^r$ we may write the different values of $\phi(x)$ in an ascending order $u_1 < \ldots < u_L$, $L\leq 2^r-1$.
Let $E_l(x)$ denote the set of all simplices $J \in \mathcal{F}^r$ (resp. $x_J$) such that $\phi[J](x) = u_l$. Then $E_1(x),\ldots, E_L(x)$ form a partition of $\mathcal{F}_r$ which we denote by $\cA(x)$. The function $x\mapsto \cA(x)$ is locally constant for Lebesgue-almost every $x\in M^r$ by a result of Chazal and Divol \cite{divol2019density} and we find 
$J_\ell \in E_\ell(x)$, $\ell=1,\ldots,L$, minimal w.r.t.\ the partial order induced by inclusion such that
$$
	u_1 = \phi[J_1](x) < u_2 = \phi[J_2](x) < \ldots < u_L = \phi[J_L](x).
$$
In particular, the functional $\rho$ from \eqref{D:RhoFunction} is Lebesgue-a.e.\ positive.
Regarding the sets $J_1,J_2,\ldots,J_L$, we remark that provided the filtration function satisfies (K1) to (K5), the existence of a minimal set is guaranteed for a.e. $x\in M^r$, see Chazal and Divol~\cite{divol2019density} Lemma 12. From now on, we just say ``minimal (w.r.t.\ inclusion)'' to characterize this set. If (K5') is satisfied instead of (K5), we have $u_1 = 0$ because all singletons appear at time 0 and we can choose $J_1 = \{1\}$ for convenience while $E_1(x) = [r]$.

Now we make the following assumption on the sublevel set of the filtration function $\phi$.
\begin{assumption}\label{A:Rho}
The inequality \eqref{E:Rho-Assumption} is satisfied for some $\alpha \in (0,1]$, viz., $\lambda^{d r}( \{ x\in M^r \mid \rho(x) \leq t\} ) \le C^* t^{\alpha}$, $t\in[0,1]$. The constant $C^*$ only depends on $r, d$ and $M$.
\end{assumption}

\begin{remark}\label{R:Rho}
Assumption~\ref{A:Rho} can be interpreted as follows: If $\rho(x) \le t$ for some $t\in (0,1]$, then we find for each $\epsilon>0$ (at least two) index sets $J_1\neq J_2$ such that
\[
	\phi[J_1](x) < \phi[J_2](x) \text{ and } \phi[J_1](y) > \phi[J_2](y)
\]
or
\[
	\phi[J_1](x) = \phi[J_2](x) \text{ and } \phi[J_1](y) > \phi[J_2](y)
\]
for some $y\in B^{\otimes r}(x;t+\epsilon)$. Hence, for $\epsilon \to 0$, we have by the continuity of the filtration function (i.e., by (K5)) either
\begin{align}\label{E:Int2}
	\phi[J_1](x) < \phi[J_2](x) \text{ and } \phi[J_1](y) \ge \phi[J_2](y)
\end{align}
or
\begin{align}\label{E:Int3}
	\phi[J_1](x) = \phi[J_2](x) \text{ and } \phi[J_1](y) > \phi[J_2](y)
\end{align}
for some $y\in B^{\otimes r}(x;t)$. It is important to note that we can assume without loss of generality that $J_2$ is minimal with respect to inclusion both in \eqref{E:Int2} and \eqref{E:Int3}. Indeed, if $J'\subseteq J_2$ is such that $\phi[J_2](x) = \phi[J'](x)$, then we can choose $J_2=J'$ right-away as $\phi[J'](y)$ still satisfies the condition on $\phi[J_1](y)$ in both cases.
In particular, we can assume without loss of generality that for any strict subset $\wt J$ of $J_2$ the condition on $x$  both in \eqref{E:Int2} and \eqref{E:Int3} is violated, viz., 
$\phi[\wt J](x) < \phi[J_1](x)$. This will be of importance later.

Note that provided $J_1\cap J_2=\emptyset$ and $J_1$ has at least two elements, we have that $\{x\in M^r: \phi[J_1](x) = \phi[J_2](x) \}$ has measure zero: Indeed by property (K5) resp. (K5') the (partial) gradient of $\phi$ is non-vanishing a.s.e., thus, $\lambda^{d|J_1|}(\{x_{J_1}: \phi[J_1](x) = a \}) = 0$ for all $0<a<\infty$. Hence, \eqref{E:Int3} can only occur if $J_1\cap J_2 \neq\emptyset$.
\end{remark}

Next, we illustrate Assumption~\ref{A:Rho} for the Vietoris-Rips filtration function.
\begin{remark}[Stability for the Vietoris-Rips filtration function]\label{R:VR}
If the Vietoris-Rips filtration is applied, we can derive a lower bound on $\rho$ by elementary means. Indeed, assume on the one hand that all pairwise distances $\|x_i-x_j\|$ are separated by more than $4t$, then the ordering of all $\|y_i-y_j\|$ does not change and all remain distinct when passing from $x$ to $y\in B^{\otimes r}(x;t)$. Thus, $\rho(x) > t$ in this case.

On the other hand, if $\rho(x) > t$, then the ordering of all $\|y_i - y_j\|$ is the same for $y \in B^{\otimes r}(x;t)$. In particular, as the filtration time $\phi[J](y)$ is determined by pairwise distances, the ordering of the entire filtration function is the same, too; viz., 
\[
	L^x = L^y \text{ and } J_i^x = J_i^y \text{ for } 1\le i \le L^x
\]
for all $y\in B^{\otimes r}(x;t)$ and for certain minimal subsets which satisfy for this filtration $J^x_1\subset J^x_2\subset \ldots \subset J^x_{L^x} \subset [r]$. Indeed, these subsets can be characterized via the pairwise distances: Assume that $x_1,\ldots,x_r$ are in general position such that all pairwise distances are distinct and $(i_1,j_1),\ldots,(i_{L-1},j_{L-1})$ are pairs such that $\|x_{i_k} - x_{j_k} \| < \|x_{i_{k+1}}-x_{j_{k+1}}\|$. Then $J^x_1 = \{i_1\}$, $J^x_2 = \{i_1,j_1\}$ and $J^x_{k+1} = \{i_{k},j_{k}\}$ for $2\le k\le L^x-1$. If some pairwise distances agree for certain pairs $(i,j)$, then we add all these pairs in a single step; plainly, the $x\in M^r$ allowing for such an equality of pairwise distances have Lebesgue-measure zero.

Conversely, assume that $\rho(x) \le t$ and $x_1,\ldots,x_r$ are in general position such that all pairwise distances are distinct and ordered via $\|x_{i_k}-x_{j_k}\| < \|x_{i_{k+1}}-x_{j_{k+1}}\|$. Then we find for each $\epsilon> 0$ a minimal index $k$ and a $y \in B^{\otimes r}(x;t+\epsilon)$ which satisfy
$$
	\|x_{i_k} - x_{j_k} \| < \|x_{i_{k+1}}-x_{j_{k+1}}\| \text{ and } \|y_{i_k} - y_{j_k} \| > \|y_{i_{k+1}}-y_{j_{k+1}}\|.
$$
This implies by the continuity of the Euclidean distance that for a minimal index $k$ and a 
$y \in B^{\otimes r}(x;t)$
$$
	\|x_{i_k} - x_{j_k} \| < \|x_{i_{k+1}}-x_{j_{k+1}}\| \text{ and } \|y_{i_k} - y_{j_k} \| \ge \|y_{i_{k+1}}-y_{j_{k+1}}\|.
$$
If $\rho(x) = t$, then we have actually equality between $\|y_{i_k} - y_{j_k} \|$ and $\|y_{i_{k+1}}-y_{j_{k+1}}\|$ in this last situation.

In total, as pairwise distances agree with Lebesgue-measure 0, the symmetric difference
\begin{align}\begin{split}\label{E:StabRhoVR}
	 &\{ x\in M^r: \rho(x)\le t \} \triangle \Big\{ x\in M^r: \exists y\in B^{\otimes r}(x;t) \text{ such that } \\
	 &\qquad \|x_i-x_j\| < \|x_k-x_\ell\| \text{ and }  \|y_i - y_j \| \ge \|y_k-y_\ell\| \text{ for some } i,j,k,\ell \Big\}
	 \end{split}
\end{align}
has Lebesgue-measure 0.
\end{remark}

The insights of the last Remark~\ref{R:VR} immediately imply
\begin{theorem}\label{L:A1forVR} Consider the Vietoris-Rips filtration function.
There is a constant $C_{VR}$ which only depends on $\diam M$ and $d$ such that for all $t\le 1$ and all $J_1 = \{i,j\}\neq \{k,\ell\} = J_2$
\begin{align}\label{E:A1forVR0}
	\lambda^{dr}(x\in M^r: \text{ \eqref{E:Int2} is valid for some } y\in B^{\otimes r}(x;t) ) \le C_{VR} \, t.
\end{align}
In particular, Assumption~\ref{A:Rho} is satisfied with 
$\alpha = 1$ and we have
\begin{align}\label{E:A1forVR1}
\lambda^{d r} \big( \{x\in M^r: \rho(x) \le t \} \big) \lesssim r^4 t.
\end{align}
\end{theorem}
\begin{proof}
We rely on characterization of $\{\rho \le t\}$ given in \eqref{E:StabRhoVR}.
Consider for two pairs $(i,j)\neq (k,\ell)$ the condition 
\[
\|x_i-x_j\| < \|x_k-x_\ell\| \text{ and }  \|y_i - y_j \| \ge \|y_k-y_\ell\| 
\]
for some $x\in M^r$ and $y\in B^{\otimes r}(x;t)$. W.l.o.g.\ assume $i\neq k\neq j$. Then
\begin{align}\label{E:A1forVR1A}
	x_k \in B(x_\ell, \|x_i-x_j\|+4t )\setminus B(x_\ell, \|x_i-x_j\|)
\end{align}
is necessarily true. Clearly, using $(b-a)^d \leq db^{d-1}(b-a)$ for $0<a<b$,
\begin{align}
	&\lambda^d\Big(B(x_\ell, \|x_i-x_j\|+4t )\setminus B(x_\ell, \|x_i-x_j\|) \Big) \label{E:A1forVR1B} \\
	&= \frac{ \pi^{d/2} }{\Gamma(d/2+1)} \cdot \Big( (\|x_i-x_j\|+4t )^d - \|x_i-x_j\|^d \Big) \nonumber \\
	&\le \frac{ \pi^{d/2} }{\Gamma(d/2+1)} \cdot d \cdot  (\diam M + 4)^{d-1} \cdot 4t. \nonumber
\end{align}
This shows \eqref{E:A1forVR0} for $C_{VR} \coloneqq \frac{ \pi^{d/2} }{\Gamma(d/2+1)} \cdot 4d \cdot (\diam M + 4)^{d-1}$. Moreover, the number of distinct pairs $(i,j)\neq (k,\ell)$ is $(r^4-2r^3-r^2+2r)/4 $. Thus, \eqref{E:A1forVR1} is satisfied if we take into account the result from \eqref{E:StabRhoVR}.
\end{proof}

Next, we study the exponent $\alpha$ for the \v Cech filtration function.
\begin{remark}
We first note that in this case the sets $J_1$ and $J_2$ in \eqref{E:Int2} and \eqref{E:Int3} are not necessarily given both by pairs. Indeed, consider $d=2$, $r=3$ and $x=\{x_1,x_2,x_3\}\subseteq M\subset \R^2$. Let $S$ be the circumsphere of $x_1,x_2$. Choose $x_3$ outside of $S$ such that $\|x_3 - x_1\| < \|x_3-x_2\| < \|x_1-x_2\|$.
Then we have
\[
	\phi[\{1,3\}](x) < \phi[\{2,3\}](x) < \phi[\{1,2\}](x) < \phi[\{1,2,3\}](x).
\]
For $y=\{y_1,y_2,y_3\}$, we set $y_1=x_1$, $y_2=x_2$ and choose $y_3$ inside the circumsphere $S$. Then
\[
	\phi[\{1,3\}](y) < \phi[\{2,3\}](y) < \phi[\{1,2\}](y) = \phi[\{1,2,3\}](y).
\]
Hence, the condition from \eqref{E:Int2} is in force for the sets $J_1 = \{1,2\}$ and $J_2=\{1,2,3\}$. Plainly, exchanging the roles of $x$ and $y$ yields the condition from \eqref{E:Int3} with the same sets. It is apparent that we can find for each $t>0$ two points $x_3,y_3$ which satisfy $\|x_3-y_3\|\le t$ additionally to the above construction.
\end{remark}

\begin{theorem}\label{L:A1forCech-Filtration}
Consider the \v Cech filtration function. There is a constant $C_{\cC}$ that only depends on $\diam M$ and $d$ such that for all $t \le 1$ and all $J_1\neq J_2$
\[
	\lambda^{dr}(x\in M^r: \text{ \eqref{E:Int2} or \eqref{E:Int3} are valid for some } y\in B^{\otimes r}(x;t) ) \le C_{\cC} t.
\]
In particular, Assumption~\ref{A:Rho} is satisfied with 
$\alpha = 1$ and we have
$$
\lambda^{d r} \big( \{x\in M^r: \rho(x) \le t \} \big) \lesssim r^{2(d+1)}t , \quad t\le 1.
$$
\end{theorem}

Consequently, the \v Cech filtration enjoys the same order of stability as the Vietoris-Rips filtration. This result seems intuitive because the filtration time in the \v Cech complex tends to vary slower than that in the Vietoris-Rips complex as the former is determined by multifold intersections while the latter is simply specified in terms of pairwise distances.

Furthermore, we see in simple examples that we cannot expect the exponent $\alpha$ to be greater than 1 for both filtration functions. Consider simply the index sets $J_1 = \{i,j\} \neq \{k,\ell\} = J_2$. Since in this case the arguments for the Vietoris-Rips and the \v Cech filtration are the same and we can rely on calculations as in \eqref{E:A1forVR1A} and \eqref{E:A1forVR1B} to see that for $x=(x_i,x_j,x_k,x_\ell)$ and $y=(y_i,y_j,y_k,y_\ell)$
\begin{align*}
	&\lambda^{4d}( \{ x\in M^r : \|x_i-x_j\| < \|x_k-x_\ell\| \text{ and } \\
	&\qquad \|y_i-y_j\| < \|y_k-y_\ell\|  \text{ for some } y \in B^{\otimes r}(x;t) \} ) \asymp t.
\end{align*}
This then shows that 
\[
	\lambda^{d r} \big( \{x\in M^r: \rho(x) \le t \} \big) \asymp t , \quad t\le 1
\]
for the Vietoris-Rips and the \v Cech filtration because a violation of the condition in~\eqref{E:Int2} with two one-dimensional simplices yields a lower bound on the exponent.\\[5pt]

A direct verification of Assumption~\ref{A:Rho} might be difficult for general filtration functions $\phi$. A sufficient condition can be established by a first order approximation. To this end, we shortly study the gradient of the Vietoris-Rips and \v Cech filtration function. Set
\begin{align}\label{E:NormGradient}
	c^* \coloneqq \max_{J\in [r]} \underset{x\in M^r}{\operatorname{ess \ sup}} \ \|\nabla \phi[J](x)\| < \infty,
\end{align}
where $\nabla \phi[J]$ denotes the gradient of $\phi[J]$.
\begin{remark}[Gradient of the Vietoris-Rips filtration function]\label{R:GradientVR}
Let $J \subseteq [r]$ and let $A$ be the set of pairs $(i,j)$ such that $i,j\in J$ and $\phi[J](x) = \|x_i - x_j\|$. Consequently, for all $(k,\ell)\in J\setminus A$, we have $\phi[J](x) > \|x_k-x_\ell\|$. Then we may identify $\nabla \phi[J](x)$ with the vector
$$
	(\nabla^1 \phi[J](x),\ldots, \nabla^r\phi[J](x)),
$$
where $\nabla^\ell \phi[J](x)$ denotes the partial gradient with respect 
to $x_\ell$.  Then we have for a pair $(i,j)\in A$ the equalities
$$	
	\nabla^i \phi[J](x)  = (x_i-x_j)/\|x_i-x_j\| = -\nabla^j \phi[J](x)
$$
while $\nabla^\ell \phi[J](x) = 0$ for $\ell$ not occuring in a pair in $A$.
Therefore and since for almost all $x\in M^r$ the set $A$ consists of a single pair $(i,j)$, it follows that
\begin{align*}
	\|\nabla \phi[J](x) \| 
	= \left( \sum_{\ell=1}^{r} \| \nabla^\ell \phi[J](x) \|^2 \right)^{1/2} 
	= \left( 1+1\right)^{1/2} 
	= \sqrt{2} = c^*.
\end{align*}
\end{remark}

\begin{remark}[Gradient of the \v Cech filtration function]\label{R:GradientCech}
Set $r=\phi[J](x)$ and let $y\in B^{\otimes r}(x;t)$. Then, on the one hand $\bigcap_{j\in J} B(y_j,r+t) \supseteq \bigcap_{j\in J} B(x_j,r) \neq \emptyset$ and 
on the other hand $\bigcap_{j\in J} B(y_j,r-t-\epsilon) \subseteq \bigcap_{j\in J} B(x_j,r-\epsilon) = \emptyset$ for all $\epsilon>0$, cf.\ Lemma \ref{L:Lemma0}.
Thus,
\begin{align}\label{E:GradientCech1}
	\phi[J](y) \in [r-t,r+t].
\end{align}
Now, let $\nu\in\R^{dr}$ with $\|\nu\|=1$ be a unit vector. Then if $\phi[J]$ is differentiable in $x\in M^r$, we have
\[
	\lim_{h\to 0} h^{-1} ( \phi[J](x+h\nu) - \phi[J](x) ) = \langle \nabla \phi[J](x), \nu \rangle.
\]
Consequently, $\| \langle \nabla \phi[J](x), \nu \rangle \| \le 1$ in each direction $\nu$. In particular, choosing $\nu = \nabla \phi[J](x) / \|\nabla \phi[J](x)\|$, we see that $c^* \le 1$.
\end{remark}

Next, we present the announced result for general filtration functions.
\begin{theorem}\label{T:GradientCondition}
Let $c^*$ be defined as in \eqref{E:NormGradient}. Then Assumption~\ref{A:Rho} is satisfied if there is an $\alpha\in (0,1]$ and a $c_0\in\R_+$ such that for any $J_1\neq J_2 $ which are no subsets of each other
\begin{align}\label{E:GradientCondition1}
\lambda^{dr}\left(	\left\{	x\in M^r:  \phi[J_2](x) \in [ \phi[J_1] (x), \phi[J_1] (x)  + 2c^* \sqrt{r}t ] \right\} \right) \le c_0 t^{\alpha}
\end{align}
and if $J_2\subset J_1$, 
\begin{align}\label{E:GradientCondition2}
	\lambda^{dr}\left( \{ x\in M^r: \phi[J_1](x) \in (\phi[J_2](x) , \phi[J_2](x) +2c^*\sqrt{r}t] \}^{(t)} \right) \le c_0 t^{\alpha}.
\end{align}
Denote $m_r = \sum_{k_1,k_2=1}^r\sum_{j=0}^{k_1\wedge k_2} \binom{r}{k_1} \binom{r-k_1}{k_2-j}$ the number of distinct subsets $J_1 \neq J_2$, then
\begin{align*}
	\lambda^{d r} \big( \{x\in M^r: \rho(x) \le t \} \big) \le c_0 m_r t^\alpha.
\end{align*}
Moreover, if the \v Cech (resp.\ Vietoris-Rips) filtration is applied, we can tighten the intervals in \eqref{E:GradientCondition1} and \eqref{E:GradientCondition2} and replace $2c^*\sqrt{r}$ by 2 (resp. by $2\sqrt{2}$).
\end{theorem} 
\begin{remark}
The number $m_r$ of distinct subsets $J_1\neq J_2$ grows exponentially in $r$. However, this is simply due to the fact that for general filtration functions we cannot derive a sharper upper bound. In case of the Vietoris-Rips and \v Cech filtration, we obtain polynomial bounds instead; we refer to the results above.

Observe that under the requirement $J_2\subset J_1$, the gradient condition in \eqref{E:GradientCondition1} just reads $\phi[J_2](x) = \phi[J_1](x)$. Consequently, we need to establish the different criterion in \eqref{E:GradientCondition2} owing to the fact that in general the measure of $\phi[J_2](x) = \phi[J_1](x)$ does not scale with $t$. Consider for example the \v Cech filtration in this setting: Given a realization of the elements $x_{J_2}$, take a realization of $x_{J_1\setminus J_2}$, whose elements all lie in the interior of the circumsphere of $x_{J_2}$. Clearly, this happens with a positive probability that is independent of $t$ and it is apparent that the filtration times of $x_{J_1}$ and $x_{J_2}$ agree in this case.
\end{remark}

\section{$L^p$-$m$-approximable feature vectors}\label{Section_FeatureVectors}

We observe a time series of point clouds $(\cX_t)_t$ in $M^r$ each of which contains exactly $r$ elements,
$$
	\cX_t = ( X_{t,1}, \ldots, X_{t,r} )^T \subseteq M^r.
$$
We build a time series of persistence diagrams $(\pd_{k,t})_t$ on top of $(\cX_t)_t$ through the filtering functions $\Filt(\cX_t)$ for a certain feature dimension $k$.

Given the point clouds $\cX_t$ of $r$ elements, let us consider a fixed feature dimension $k \in \{0,1,\ldots,d-1\}$ such that $k+2 \le r$. (A sufficient condition is clearly that $r\ge d+1$.) 

Then the random number of $k$-dimensional features in $\pd_{k,t}$ admits the deterministic upper bound $N_k \coloneqq \binom{r}{k+1} $ (see also the textbook of Boissonnat et al.~\cite{boissonnat2018geometric}), i.e.,
\begin{align*}
	&\#\{ \text{$k$-dimensional features in }\pd_{k,t} \} \le N_k  \le r^{k+1}.
\end{align*}
(Plainly, a refined bound is given by $r(r-1)$ for the Vietoris-Rips filtration.) 

W.l.o.g.\ we can assume that the actual number of features in $\pd_{k,t}$ equals $N_k$ by borrowing trivial features from the diagonal $\Delta = \{(x,x): x\in [0,\infty) \}$.
Now, we study feature vectors obtained from the corresponding $k$th persistence diagrams $(\pd_{k,t})_t$ of $\cX_t$: Let 
$$
(b_{t,1},d_{t,1}), \ldots, (b_{t,N_k},d_{t,N_k})
$$
be the features of $\pd_{k,t}$ given in ascending order according to their birth, i.e.,
$$
	b_{t,1} \le b_{t,2} \le \ldots \le b_{t,N_k}.
$$
(If there are features with the same birth, we order these at random.) Then, we study the following representation of $\pd_{k,t}$
\begin{align}\label{E:RepresentationPD1}
	Z_{k,t} = 
	\begin{pmatrix}
		b_{t,1} \\
		d_{t,1} \\
		\vdots \\
		b_{t,N_k}\\ 
		d_{t,N_k} \\
	\end{pmatrix}.
\end{align}
In order to focus on certain topological aspects we consider ``aggregated'' statistics $f(Z_{k,t})$, where $f\colon \R^{2N_k} \to \R^{\ell}$
is some Lipschitz-continuous function.
 
\begin{example}[Total $\gamma$-persistence of the persistence diamgram]
Given $Z_{k,t}$ and $\gamma\in [1,\infty]$ the total $\gamma$-persistence is defined by
$$
	\| (d_{t,i}-b_{t,i})_{i=1,\ldots,N_k} \|_{\gamma} = \left(\sum_{i=1}^{N_k} (d_{t,i}-b_{t,i})^{\gamma}\right)^{1/\gamma}
$$	
resp., 
$$
	\| (d_{t,i}-b_{t,i})_{i=1,\ldots,N_k} \|_{\infty} = \max_{i\le N_k} ( d_{t,i}-b_{t,i} ).
$$
\end{example}
Next, we give more details on the data generating procedure. The observations are assumed to be stationary and weakly dependent in the following sense (and this will be the null hypothesis in the next section). 
\begin{assumption}[$L^p$-$m$-approximable data]\label{Assumption:Lrm-Condition-PointProcess}
The point cloud process $\cX=(\cX_t)_t$ has a Bernoulli shift structure of the form 
\begin{align}\label{Eq:BernoulliShiftStructure}
	\cX_t = ( X_{t,1}, \ldots, X_{t,r} ) = F(\epsilon_t,\epsilon_{t-1},\epsilon_{t-2},\ldots)	\qquad (t\in\Z),
\end{align}
where $(\epsilon_t)_t$ are iid elements taking values in some measurable space $S$, and $F=(F_0,\ldots,F_{r-1})\colon S^{\infty} \to M^r$ is a measurable function. 
Moreover, we assume $\cX_t$ admits a bounded density on $M^r$.
Let $(\cmX_t)_t$, resp. $(X_{t,i}^{(m)})_t$, denote the $m$-dependent process constructed from
\begin{align}\label{Eq:m-dep-Approx}
	X_{t,i}^{(m)} \coloneqq F_i(\epsilon_t,\epsilon_{t-1},\ldots,\epsilon_{t-m+1},\epsilon_{t-m}^{(t)},\epsilon_{t-m-1}^{(t)},\ldots ) \qquad (t\in\Z),
\end{align}
where here and in the following for each $t\in\Z$ the sequence $(\epsilon_s^{(t)})_s$ is an independent copy of the innovations $(\epsilon_s)_s$.
Granted $\sum_{i=1}^r \E[ \|X_{0,i}\|^p ] < \infty$ for a certain $p\ge 1$ we assume $(\cX_t)_t$ satisfies the $L^p$-$m$-approximation property 
\begin{align}\label{Eq:L^p-m-property}
	\sum_{m\ge 1} \gamma_p(m)^{\alpha/(\alpha+p)} < \infty,
\end{align} 
where $\gamma_p(m) \coloneqq \sum_{i=1}^r \E[\| X_{0,i} - X_{0,i}^{(m)}\|^p] ^{1/p}$
and $\alpha>0$ is determined by the filtration as described in Assumption~\ref{E:Rho-Assumption}.
This completes the assumption.\\[5pt]
\end{assumption}
We give an example for data structures which satisfy this assumption.
\begin{example}[Point clouds obtained from time series via delay embedding]\label{E:Bernoulli}
Suppose $Y=(Y_t)_t$ is a stationary multivariate $M$-valued time series of the type
\begin{align*}
	Y_t = G(\epsilon_t,\epsilon_{t-1},\ldots),
\end{align*}
where $(\epsilon_t)_t$ is a sequence of iid elements in some measurable space $S$ and $G\colon S^{\infty} \to M$ is a measurable function. 
Let us assume $\E[ \|Y_0\|^p] < \infty$ and
\begin{align}\label{Eq:Lrm-ConditionForTS}
	M_0 \coloneqq \sum_{m=0}^{\infty}  \ \E[ \| Y_{0} - G(\epsilon_{0},\ldots, \epsilon_{-m+1}, \epsilon_{-m}^{(0)}, \epsilon_{-m-1}^{(0)},\ldots ) \|^p ]^{\frac{\alpha}{p(\alpha+p)}} < \infty,
\end{align}
where the exponent $\alpha$ is again determined through the filtration by \eqref{E:Rho-Assumption}.

We can then consider an associated point cloud process 
\begin{align*}
	\cX_t &\coloneqq (X_{t,1},\ldots,X_{t,r}) \coloneqq (Y_t,Y_{t-1},\ldots,Y_{t-r+1}),
\end{align*}
i.e.,
\[
	X_{t,i} = Y_{t-i+1} = G(\epsilon_{t-i+1},\epsilon_{t-i},\ldots).
\]
In this case $F_i(a_t,a_{t-1},\ldots) = G(a_{t-i+1},a_{t-i},\ldots)$ and the $m$-dependent approximation \eqref{Eq:m-dep-Approx} is given by
\begin{align}
	X_{t,i}^{(m)} 
	&= F_i(\epsilon_t,\ldots,\epsilon_{t-m+1},\epsilon_{t-m}^{(t)},\epsilon_{t-m-1}^{(t)},\ldots ) \nonumber\\
	&=
	\begin{cases}\label {E:Ex1} 
		G( \epsilon_{t-i+1},\ldots, \epsilon_{t-m+1}, \epsilon^{(t)}_{t-m},\ldots )  & \text{ if } i \leq m \\
		G(\epsilon^{(t)}_{t-i}, \epsilon^{(t)}_{t-i-1}, \ldots ) & \text{ if } i > m.
	\end{cases}
\end{align}
Recall that $\gamma_p(m) =   \sum_{i=1}^r \E[ \|X_{0,i}-X_{0,i}^{(m)}\|^p ]^{1/p}$. Thus, if $m\ge r$, we use that $X_{0,i}=Y_{-i+1}$ and infer from \eqref{E:Ex1} that
\begin{align*}
	X_{0,i}-X_{0,i}^{(m)}  =  G(\epsilon_{-i+1},\epsilon_{-i},\ldots) - G( \epsilon_{-i+1},\ldots, \epsilon_{-m+1}, \epsilon^{(0)}_{-m},\ldots ) .
\end{align*}
Hence, due to \eqref{Eq:Lrm-ConditionForTS} and stationarity,
\begin{align*}
	&\sum_{m=r}^\infty \gamma_p(m)^{\alpha/(\alpha+p)} 
	\le \sum_{m=r}^\infty  \sum_{i=1}^r \E[ \|X_{0,i}-X_{0,i}^{(m)}\|^p ]^{\frac{\alpha}{p(\alpha+p)}}\\
	&=   \sum_{m=r}^\infty \sum_{i=1}^r \E[ \| G(\epsilon_{-i+1},\epsilon_{-i},\ldots) - G( \epsilon_{-i+1},\ldots, \epsilon_{-m+1}, \epsilon^{(0)}_{-m},\ldots )  \|^p ]^{\frac{\alpha}{p(\alpha+p)}} \\
	 \\
	&=    \sum_{i=1}^r \sum_{m=r}^\infty \E[ \| Y_{0} - Y_{0}^{(m-i+1)}  \|^p ]^{\frac{\alpha}{p(\alpha+p)}} 
	 \le r M_0 < \infty,
\end{align*}
because $m\ge r\ge i$. In particular, $\sum_{m\ge 1} \gamma_p(m)^{\alpha/{\alpha+p}}$ is finite.
To illustrate, let us consider a classical linear process 
$$
	Y_t = \sum_{k=0}^{\infty} a_k\epsilon_{t-k} = G(\epsilon_t,\epsilon_{t-1},\ldots)
$$
with $a_k\in\R^{d\times d}$. Then, $X_{t,i} = Y_{t-i} = \sum_{k=0}^\infty a_k\epsilon_{t-i-k}$ and according to \eqref{E:Ex1},
\begin{align*}
	X_{t,i}^{(m)} 
	&=
	\begin{cases}
		\sum_{k=0}^{m-i} a_k \epsilon_{t-i+1-k} + \sum_{k=m-i+1}^{\infty} a_k\epsilon_{t-i+1-k}^{(t)}  & \text{,if } i \leq m \\
		\sum_{k=0}^{\infty} a_k \epsilon_{t-i-k}^{(t)} & \text{,if } i > m.
	\end{cases}
\end{align*}
Furthermore, 
\begin{align*}
	M_0 
	&= \sum_{m=1}^\infty \Big(\E[ \| \sum_{k\geq m} a_k(\epsilon_{-k} - \epsilon_{-k}^{(0)}) \|^p ]^{1/p}\Big)^{\frac{\alpha}{\alpha + p}}\\
	&\leq \sum_{m=1}^\infty \Big( 2\E[\|\epsilon_0\|^p]^{1/p} \sum_{k\geq m} \|a_k\| \Big)^{\frac{\alpha}{\alpha + p}}\\
	&\leq (2\E[\|\epsilon_0\|^p]^{1/p})^{\frac{\alpha}{\alpha + p}}\, \sum_{k\geq 1} k \|a_k\|^{\frac{\alpha}{\alpha + p}},
\end{align*}
so $\sum_{k\geq 1} k \|a_k\|^{\frac{\alpha}{\alpha + p}} < \infty$ is sufficient to ensure \eqref{Eq:Lrm-ConditionForTS} and, consequently, \eqref{Eq:L^p-m-property}.
\end{example}

Next, we present the main probabilistic results of the manuscript regarding the weak dependence of the feature vectors $(Z_{k,t})$. These are all formulated for point cloud data which satisfies Assumption~\ref{Assumption:Lrm-Condition-PointProcess} which itself relies on Assumption~\ref{A:Rho}.

\begin{theorem}\label{T:LpmApprox}
Let $k\le r+2$ and suppose $(\cX_t)_t$ satisfies Assumption \ref{Assumption:Lrm-Condition-PointProcess} for some $p\geq 1$ and for a filtering function which satisfies Assumption~\ref{A:Rho}. Then the corresponding process of feature representations $(Z_{k,j})_j$ is $L^p$-$m$-approximable in the sense that 
\begin{align*}
	\sum_{m=1}^{\infty} \E [\| Z_{k,0} - Z^{(m)}_{k,0} \|^p] ^{1/p} < \infty.
\end{align*}
\end{theorem}

The last Theorem~\ref{T:LpmApprox} enables us to obtain a functional central limit theorem for data derived from the feature representations $(Z_{k,t})_t$. To this end, denote $(W_\Gamma (t): t\in [0,1])$ an $\ell$-dimensional Brownian motion with covariance matrix $\Gamma\in\R^{\ell\times\ell}$, i.e., a Gaussian process with mean $0\in\R^\ell$ and covariance function $Cov( W_\Gamma(s), W_\Gamma(t) )= \min\{s,t\} \Gamma$.
 
\begin{theorem}\label{Theorem:FCLT}
Let the requirements of the previous Theorem~\ref{T:LpmApprox} be in force for $p\ge 2$. Let $\ell\in\N$ and let 
$f\colon \R^{2N_k} \to \R^{\ell}$ be a Lipschitz-continuous function. Then
the series $\Gamma_k = \sum_h \Cov( f(Z_{k,h}),f(Z_{k,0}) )$ converges (coordinatewise) absolutely and 
\begin{align*}
		\frac{1}{\sqrt{n}} \sum_{j=1}^{[n \cdot ]} \Big( f(Z_{k,j}) - \E f(Z_{k,0}) \Big) \Rightarrow W_{\Gamma_k} \quad (n\to\infty),
\end{align*}
where the weak convergence takes place in the $\ell$-dimensional Skorohod space $D_{\ell}([0,1])$.
\end{theorem}	

In order to detect a structural break, we will rely on the following result.
\begin{cor}\label{Cor:FCLT}
Assume the conditions of Theorem~\ref{Theorem:FCLT}. Then
\begin{align*}
	\frac{1}{\sqrt{n}} \left( \sum_{j=1}^{[n \cdot ]} f(Z_{k,j}) - \ \cdot\times \sum_{j=1}^n f(Z_{k,j}) \right) \Rightarrow W_{\Gamma_k} - \ \cdot \times W_{\Gamma_k}(1) =: B_{\Gamma_k} \quad (n\to\infty),
\end{align*}
weakly in the $\ell$-dimensional Skorohod space $D_{\ell}([0,1])$.
\end{cor}

\section{Testing for structural breaks in a sequence of PDs}\label{Section_Test}
We test the null hypothesis that the persistence diagrams have a constant distribution, i.e., 
\begin{align}\label{E:H0}
	H_0\colon \cL(f(Z_{k,1})) = \ldots = \cL(f(Z_{k,n})),
\end{align}
against the alternative ``exactly one structural break occurs'', i.e., 
\begin{align}
\begin{split}\label{E:H1}
	H_1\colon &\cL(f(Z_{k,1})) = \ldots = \cL(f(Z_{k,v^*})) \\
	&\not= \cL(f(Z_{k,v^*+1})) = \ldots = \cL(f(Z_{k,n})).
\end{split}\end{align}
Here $v^*$ denotes the (typically unknown) time index where the break occurs. To construct a test for distinguishing between $H_0$ and $H_1$, we consider the vectors $(f(Z_{k,j}))_j$, where $f\colon\R^{2N_k}\to\R^\ell$ is Lipschitz, and use the cumulative sum (CUSUM) statistic 
\begin{align*}
	S_v = \frac{1}{\sqrt{n}} \left( \sum_{j=1}^v f(Z_{k,j}) - \frac{v}{n}\sum_{j =1}^n f(Z_{k,j})  \right) \qquad (1\le v\le n),
\end{align*}
which compares the overall mean with the mean of $v$ observations. The following test statistics are established in practice (see for instance Aue et al.\ \cite{aue2009break})
\begin{align*}
	\Lambda_{k,n} \coloneqq \max_{1\leq v\leq n} S_v^T \wh{\Gamma}_{n,k}^{-1} S_v, \qquad \Omega_{k,n} \coloneqq \frac{1}{n} \sum_{v=1}^n S_v^T \wh{\Gamma}_{n,k}^{-1} S_v,
\end{align*}
where $\wh{\Gamma}_{n,k}$ is an estimator for the long run covariance matrix of the process $(f(Z_{k,j}))_j$ which is defined by $\Gamma_k = \sum_{h\in\Z} \Cov( f(Z_{k,h}),f(Z_{k,0}) )$. We assume that the estimator satisfies 
\begin{align}\label{C:CovEst}
\|\wh{\Gamma}_{n,k} - \Gamma_k \| = o_{a.s.}(1).
\end{align} 

Under $H_0$ the limit distributions follow from an existing result under weak dependence given in Aue et al.\ \cite{aue2009break}.
\begin{theorem}[Asymptotics under $H_0$; an application of Theorem 2.1 in Aue et al.\ \cite{aue2009break}]\label{T:AueBreak}
Assume the conditions of Theorem~\ref{Theorem:FCLT}. Moreover, let $\wh{\Gamma}_{n,k}$ be such that \eqref{C:CovEst} is in force. Then under the null hypothesis (formulated in \eqref{E:H0})
\begin{align*}
	\Lambda_{k,n} &\Rightarrow \sup_{ 0\leq t \leq 1 } \sum_{i=1}^{\ell} B_i^2(t) =: \Lambda(\ell) \qquad (n\to\infty),\\
	\Omega_{k,n} &\Rightarrow \sum_{i=1}^{\ell} \int_0^1  B_i^2(t) \,\diff t =: \Omega(\ell) \qquad (n\to\infty),
\end{align*} 
where $B_1,\ldots,B_\ell$ are independent standard Brownian bridges on $[0,1]$.
\end{theorem}
The size $r$ of the point clouds can be very large, clearly this implies that $N_k$ is large, too. Then in order not to lose information, it is also reasonable to choose $\ell$ large. Hence, the limiting behavior of the statistics
$\Lambda(\ell)$, $\Omega(\ell)$ for $\ell \to \infty$ is of interest.
Then as  detailed in Remark 2.1 in Aue et al.\ \cite{aue2009break} the standardized variables have a normal distribution in the limit, viz.,
\begin{align*}
	\ol{\Lambda}(\ell) = \frac{\Lambda(\ell) - \ell/4}{(\ell/8)^{1/2}} \Rightarrow  N(0,1) \text{ and } \ol{\Omega}(\ell) = \frac{\Omega(\ell) - \ell/6}{(\ell/45)^{1/2}} \Rightarrow  N(0,1) \quad (\ell\to\infty).
\end{align*}
Theorem~\ref{T:AueBreak} allows for asymptotic test procedures. We sketch this for the statistic $\Lambda_{k,n}$. Let $\alpha\in (0,1)$ and let $c_{\ell,\alpha}$ be the $(1-\alpha)$-quantile of $\Lambda(\ell)$. Clearly, the distribution of $\Lambda(\ell)$ is nonnegative and absolutely continuous w.r.t.\ the Lebesgue-measure. Thus, under the null hypothesis  $H_0$ of no structural break and by the Portemanteau-Theorem
\[
	\lim_{n\to\infty} \p_{H_0}( \Lambda_{k,n} > c_{\ell,\alpha} ) = \p_{H_0}(\Lambda(\ell) > c_{\ell,\alpha} ) = \alpha.
\] 
Consequently, the test
\[
		\text{ reject $H_0$ if $\Lambda_{k,n} > c_{\ell,\alpha}$ }
\]
has asymptotic level $\alpha$.
Observe that for large $\ell$ the value $(c_{\ell,\alpha} - \ell/4)/\sqrt{\ell/8}$ should be close to the $(1-\alpha)$-quantile of the normal distribution. For small and moderate values of $\ell$, the value $c_{\ell,\alpha}$ can be approximated with Monte Carlo methods.
The asymptotic Gaussianity can typically be leveraged to obtain a validity results for bootstrap and subsampling techniques; however we do not study this any closer here.

Following Aue et al.\ \cite{aue2009break}, let $\theta\in (0,1)$ be arbitrary but fixed. Set $v^* = \floor{n\theta}$.
We observe a pre-change sequence $(f(Z_{k,j}))_{j}$ on $1\le j\le v^*$ and a post-change sequence $(f(Z^*_{k,j}))_{j}$ on $v^* < j \le n$ built on top of the feature representations as outlined above. We assume that both sequences are stationary and $L^2$-$m$-approximable; as before a sufficient condition is formulated in Assumption~\ref{Assumption:Lrm-Condition-PointProcess}. 

\begin{theorem}[Asymptotics under $H_1$]\label{T:AsympH1}
Suppose the sequences $(f(Z_{k,j}))_{j}$ and $(f(Z^*_{k,j}))_{j}$ are stationary and $L^2$-$m$-approximable. Then it holds under the alternative \eqref{E:H1} that
\begin{align*}
	\sup_{t\in [0,1]} \Big| \frac{1}{\sqrt{n}} S_{\floor{nt}} - S^*(t) \Big| \rightarrow 0 \text{ in probability},
\end{align*}
where
\begin{align*}
	S^*(t) =
	\begin{cases}
		t(1-\theta) \E[ f(Z_{k,1}) - f(Z^*_{k,1})], \qquad t\in [0,\theta],\\
		\theta(1-t)  \E[ f(Z_{k,1}) - f(Z^*_{k,1})], \qquad t\in (\theta,1].
	\end{cases}
\end{align*}
\end{theorem}
We conclude this section with a short application of the asymptotic behavior of a change point estimator under $H_1$ given that there is a change in the mean.
The model is as follows. The series $Y_j = f(Z_{k,j}) \in \R^\ell$, $1\le j\le n$, undergoes a change in the mean according to
\begin{align}\label{E:StructuralBreak}
	Y_j =
	\begin{cases}
		e_j + \mu, &\qquad j\le v^*,\\
		e_j + \mu - \Delta , &\qquad j>v^*,
	\end{cases}
\end{align}
for certain $\Delta, \mu\in\R^{\ell}$ such that $\Delta \neq 0$.

Before we inspect potential change point estimators, we interpret this model in our context of the definition of $f(Z_{k,j})$. Recall that $Z_{k,j}$ is built on top of the filtered point cloud data $\cX_j$ which itself is defined through a Bernoulli system as in Example~\ref{E:Bernoulli}. Relying on Theorem~\ref{T:LpmApprox} and the Lipschitz-continuous map $f$, the time series $f(Z_{k,j})$ ($1\le j\le n$) is $L^p$-$m$-approximable. So, it can be treated with standard methods from time series analysis. However, a peculiarity of the present situation is that it is hard to establish a meaningful link between the abstract change in mean $\Delta$ and the corresponding structural change in the underlying point cloud data $(\cX_j)_j$. It seems very plausible that the same $\Delta$ can be the outcome of two quite different structural breaks in the underlying point cloud process. So, a deeper understanding of these topological methods is required in order to assess this relation any further; however, this is beyond the scope of this manuscript.

We assume that the zero mean sequence $(e_j)_j$ converges at rate $a_n$ to the mean $0\in\R^\ell$, viz.,
\begin{align}
	&a_n^{-1} \sum_{j=1}^n e_j \rightarrow 0 \ a.s. \quad (n\to\infty) \label{E:SLLNInno}
\end{align}
for a non-decreasing sequence $(a_n)_n$ with limit $+\infty$. We also write $ \sum_{j=1}^n e_j = o_{a.s.}(a_n)$ in this case.
The corresponding CUSUM statistic is
\begin{align*}
	 S_v &= \frac{1}{\sqrt{n}} \Big( \sum_{j=1}^v Y_i - \frac{v}{n} \sum_{i=1}^n Y_i \Big)=
	\begin{cases}
		\frac{1}{\sqrt{n}} \cE_v +\frac{v(1-\theta)}{\sqrt{n}}   \Delta, &\quad v \le v^*, \\
		\frac{1}{\sqrt{n}} \cE_v +  \frac{(n-v)\theta}{\sqrt{n}}\Delta  , &\quad v>v^*.
	\end{cases}
\end{align*}
Here
$$
	\frac{ 1}{\sqrt{n}} \cE_v = \frac{1}{\sqrt{n}} \Big(\sum_{j=1}^v e_j - \frac{v}{n} \sum_{j=1}^n e_j \Big) \quad (1\le v\le n),
$$
is the CUSUM process of the innovations $(e_j)_j$.
In this setting $S^*(t) = t(1-\theta) \Delta$ if $t\in [0,\theta]$ and $S^*(t) = (1-t)\theta \Delta$ if $ t\in(\theta,1]$.

Let $\wt\Sigma_n$ be positive definite, then the break point $\theta$ satisfies
$$
	\theta = \underset{t\in [0,1]}{\operatorname{arg \ max}} ~ S^*(t)^T \wt\Sigma_n S^*(t). 
$$
Since $\max_{v\le n} | n^{-1/2} S_v - S^*(v/n)| = o_p(1)$ by Theorem~\ref{T:AsympH1},
a reasonable estimator of the break point $\theta$ is given by
\begin{align}\label{E:whTheta}
    	\wh\theta_n = \frac{1}{n} \ \underset{1\le k\le n}{\operatorname{arg \ max}} ~ S_k^T \wt \Sigma_n S_k.
\end{align}
A typical choice is $\wt\Sigma_n \coloneqq \wh{\Gamma}_{n,k}^{-1}$ assuming that $\wh{\Gamma}_{n,k} \to \Gamma_{k}$ a.s. as above in \eqref{C:CovEst}. 
In this case, relying on Theorem~\ref{T:AueBreak} we can control the fluctuations and have under the null hypothesis
$$
	S_{n \wh\theta_n}^T \wt \Sigma_n S_{n \wh\theta_n} = \Lambda_{k,n} \Rightarrow \Lambda(\ell).
$$
Rescaling with $\wh{\Gamma}_{n,k}^{-1}$ puts approximately the same weight on each coordinate.
With these preparations, we formulate the following statement on the consistency of the break point estimator.

\begin{theorem}\label{T:HatTheta}
Assume the model~\eqref{E:StructuralBreak} with underlying rate of convergence $(a_n)_n$ for the innovations $(e_j)_j$ as detailed in \eqref{E:SLLNInno} such that $a_n = o(n)$. Moreover, suppose there is a measurable set $\Omega_0$ with probability 1 such that
$$
		\sup_{n\in\N} \frac{ \| \wt\Sigma_n(\omega) \| }{ \Delta^T \wt\Sigma_n(\omega) \Delta } < \infty \quad \text{for all } \omega\in\Omega_0.
$$	
Then the estimator $\wh\theta_n$ from \eqref{E:whTheta} satisfies $\wh\theta_n - \theta =  O_{a.s.}( n^{-1} a_n)$.
\end{theorem}
For instance, if the innovations $(e_j)_j$ are univariate and satisfy a law of the iterated logarithm any sequence $(a_n)_n$ such that $n^{1/2} \log\log n = o(a_n)$ yields a strongly consistent estimator $\wh\theta_n$. And the law of the iterated logarithm holds if the innovation process satisfies a strong invariance principle at a suitable rate, viz., if it can be coupled together with a standard Brownian motion $B$ on a common probability space in the sense that
\begin{align}\label{E:KMT}
	\sum_{j=1}^n e_j = \sigma B(n) + O_{a.s.}(r_n)
\end{align}
for a suitable remainder sequence $(r_n)_n$ and $\sigma^2 = \sum_{j\in\Z} \E[ e_0 e_j]$. The remainder is usually of the type $n^{1/p}$ or $\log n$ depending on moment conditions and dependence settings, we refer to Berkes et al. \cite{berkes2014komlos} for more details.

\section{Technical details on $L^p$-$m$-approximable feature vectors}\label{Section_Proofs}

\subsection{Stability of filtration functions}

Depending on the exact type of the filtration, e.g., \v Cech  and Vietoris-Rips filtration, we can have $\phi[J]\equiv 0$ for singletons $J=\{j\}$. In this case it suffices to consider $J \in [r]$ of size larger than one in the following. In all other cases, where $\phi$ is not zero for singletons, we can assume $J$ to be an arbitrary subset instead.
Furthermore, given the conditions (K1) to (K5), resp.\ (K1) to (K4) and (K5') for the filtration, the function $\phi[J]$ is differentiable with nonvanishing gradient almost everywhere on $M^r$, see \cite{divol2019density}. 

\begin{lemma}\label{L:Lemma0}
Consider the \v Cech filtration function. Let $t>0$, $x\in M^r$, $y\in B^{\otimes r}(x;t)$ and $J\subseteq [r]$ with $|J|\le d$.
\begin{itemize}
\item [(a)] Then $\phi[J](y)\in [\phi[J](x)-t,\phi[J](x)+t]$.
\item [(b)] If $j\notin J$ and $x_j$ is not contained inside the circumsphere of $x_J$, then $\phi[J](x) < \phi[J\cup\{j\}](x)$. 
\end{itemize}
\end{lemma}
\begin{proof}
Set $R=\phi[J](x)$ and let $z$ be the center of the circumsphere of the $x_J$. Recall that this center is uniquely determined. Moreover, $x_J\subseteq B(z,R)$ and $x_J \not\subseteq B(\wt z,\wt R)$ for any $\wt R< R$ and any $\wt z$.

(a) We have for each $s\ge t$ and $R'>0$
$$	B(y_j,R'-s) \subseteq B(x_j,R')\subseteq B(y_j,R'+s). $$
Thus, on the one hand $\{z\} = \bigcap_{j\in J} B(x_j,R) \subseteq \bigcap_{j\in J} B(y_j,R+t) $, yielding the upper bound.

On the other hand, $\bigcap_{j\in J} B(y_j,R-t-\epsilon) \subseteq \bigcap_{j\in J} B(x_j,R-\epsilon) = \emptyset$ for each $\epsilon >0$. This yields the lower bound.

(b) Assume that $R=\phi[J\cup\{j\}](x)$ is true and let $z'$ be the circumcenter of $x_{J\cup\{j\}}$. Then $x_{J}$ is contained in $B(z', R) \cap B(z,R)$. By the minimality of $R$ and the uniqueness of $z$, we must have $z'=z$. This implies in turn that $x_j\in B(z,R)$ which is a contradiction. Hence, $R<\phi[J\cup\{j\}](x)$. 
\end{proof}

\begin{proof}[Proof of Theorem~\ref{L:A1forCech-Filtration}]
If $\rho(x)\le t$ then equation~\eqref{E:Int2} or \eqref{E:Int3} are valid for certain sets $J_1\neq J_2$. We distinguish three cases depending on the possible inclusions between $J_1$ and $J_2$. We use often that $|\phi[J](x) - \phi[J](y)|\le t$ for all $y\in B^{\otimes r}(x;t)$ by Lemma~\ref{L:Lemma0}.

Before we work on these three cases, we give a simple upper bound on the number of potential subsets $J_1,J_2$. Recall that the number of $k$-dimensional simplices is $\binom{r}{k}$ and that the filtration times of the \v Cech filtration are determined by simplices of dimension $k\le d+1$.
Hence, the number of sets $J_1\neq J_2$ which satisfy \eqref{E:Int2} or \eqref{E:Int3} is at most
\begin{align*}
	\Big(\sum_{k=2}^{(d+1)\wedge r} \binom{r}{k} \Big)^2 \le r^{2(d+1)};
\end{align*}
the 0-dimensional simplices do not occur on the left-hand side as their filtration time is constant 0. Depending on the exact case, the relevant number of sets can be significantly less; this is however not important in the following.

Note that the sets $J_1,J_2$ which we study in the sequel contain at most $d+1$ elements, in particular, all constants that will occur only depend on $d$ and $M$ but not on $r$.

Subcase I: $J_1 \cap J_2 = \emptyset$. In this case only \eqref{E:Int2} is relevant because \eqref{E:Int3} happens with measure 0 if $J_1$ and $J_2$ are disjoint. We rely on
\begin{align*}
	&\{x\in M^r:  \phi[J_1](x) < \phi[J_2](x), \quad  \phi[J_1](y) \ge \phi[J_2](y) \text{ for one } y\in B^{\otimes r}(x;t) \} \\
	&\subseteq \{x\in M^r: \phi[J_2](x) \in ( \phi[J_1](x), \phi[J_1](x)+2t) \}.
\end{align*}
Moreover, if the \v Cech filtration function is applied Lemma 6.10 in Divol and Polonik \cite{divol2019choice} yields for $0<a<b<\infty$ an upper bound as follows
\[
	\lambda^{dr}( \{ x\in M^r: \phi[J_2](x) \in (a,b) \} ) \le c (b-a)
\]
for a constant $c\in\R_+$.
So by conditioning on the coordinates $x_{J_1}$, we find
\[
	\lambda^{dr}( \{x\in M^r: \phi[J_2](x) \in [\phi[J_1](x),  \phi[J_1](x)+2t] \} ) \le 2t.
\]
Subcase II: $J_2 \subset J_1$ or $J_1\subset J_2$. Remember that we denote the smallest closed bounding ball of $x_{\wt J}$ by $B^*(x_{\wt J})$ and its center by $z(x_{\wt J})$. Note that not all elements of $x_{\wt J}$ lie necessarily on the boundary of $B^*(x_{\wt J})$. 

We begin with the situation in \eqref{E:Int3} where $J_2\subset J_1$; then this condition reads
\[
	\phi[J_1](x) = \phi[J_2](x) \text{ and } \phi[J_1](y) > \phi[J_2](y).
\]
We may assume w.l.o.g.\ that $J_1 = J \cup \{ j_0\}$ while $J_2 = J$. 

Hence, it is sufficient to study the set
\[
	\{ x\in M^r \ | \  x_{j_0}\in B^*(x_{J}), y_{j_0}\notin B^*(y_{J}) \text{ for some } y\in B^{\otimes r}(x;t) \}.
\]
Then the condition $y_{j_0}\notin B^*(y_{J})$ implies that $x_{j_0}\notin B^*(y_{J})^{(-t)}$. Thus, we consider conditional on $x_{J}$
\begin{align}\label{E:Inclusion1}
	\lambda^d(\{x_{j_0}\in M \ | \ x_{j_0} \in B^*(x_{J}) \setminus B^*(y_{J})^{(-t)} \text{ for some } y\in B^{\otimes r}(x;t)\}).
\end{align}
Note that we can replace from now on $y\in B^{\otimes r}(x;t)$ by $ y_{J}\in B^{\otimes r}(x_{J};t)$

To study the measure of this set difference, we make the abbreviations
\[
	z_1 = z( x_J ), z_2 = z( y_J), R_1 = \phi[J](x), R_2 = \phi[J](y).
\]
Firstly, we use that $R_1 - t \le R_2 \le R_1+t$ and secondly that $A\setminus C\subseteq (A\setminus B) \cup (B\setminus C)$ in order to deduce
\begin{align*}
	 &B^*(x_{J}) \setminus B^*(y_{J})^{(-t)} = B(z_1,R_1)\setminus B(z_2,R_2-t)\subseteq B(z_1,R_1)\setminus B(z_2,R_1-2t)\\
	 &\subseteq \big( B(z_1,R_1)\setminus B(z_1,R_1-2t)\big) \cup \big( B(z_1,R_1-2t)\setminus B(z_2,R_1-2t)\big). 
\end{align*}
Note that only the second set on the right-hand side depends on $y_J$.

Hence, we can bound above \eqref{E:Inclusion1} with the sum of the following two terms
\begin{align}\label{E:Inclusion1b}
	\lambda^d(  B(z_1,R_1)\setminus B(z_1,R_1-2t) ) \le \frac{ \pi^{d/2} }{\Gamma(d/2+1)} \cdot d \cdot  (\diam M/2)^{d-1} \cdot 2t
\end{align}
and 
\begin{align}\label{E:Inclusion2}
	\lambda^d\big( \{ x_{j_0}\in  B(z_1,R_1-2t)\setminus B(z_2,R_1-2t) \text{ for some } y_J \in B^{\otimes r}(x_J;t)\} \big),
\end{align}
where $z_2 = z(y_J)$ depends on $y_J$. So, it remains to study \eqref{E:Inclusion2}. Clearly, we have by translation invariance of the Lebesgue-measure and by an obvious subset relation that \eqref{E:Inclusion2} is at most
\begin{align}
	& \lambda^d\big( \{ x_{j_0}\in  B(z_2-z_1,R_1-2t)\setminus B(0,R_1-2t) \text{ for some } y_J \in B^{\otimes r}(x_J;t)\} \big) \nonumber \\
	&\le\lambda^d\big( B\big(0,R_1-2t+\sup_{y_J \in B^{\otimes r}(x_J;t)} \|z_2-z_1\| \big)\setminus B(0,R_1-2t) \big) \nonumber\\
	&\le \frac{ \pi^{d/2} }{\Gamma(d/2+1)} \cdot d \cdot  (\diam M/2)^{d-1} \cdot \sup_{y_J \in B^{\otimes r}(x_J;t)} \|z_2-z_1\|. \label{E:Inclusion3}
\end{align}

Before taking into consideration the $\sup$ in \eqref{E:Inclusion3}, we catch up with the details of the situation given in \eqref{E:Int2} which reads
\begin{align*}
	\phi[J_1](x) < \phi[J_2](x) \text{ and } \phi[J_1](y) = \phi[J_2](y)
\end{align*}
Since in this case the roles of $J_1, J_2$ and $x,y$ are just reversed, we can assume w.l.o.g.\ $J=J_1\subset J_2=J\cup\{j_0\}$.
Hence, we arrive at
\[
	x_{j_0} \in B^*(y_{J_1})^{(t)} \setminus B^*(x_{J_1});
\]
this implies in turn
\[
	 x_{j_0} \in \big(B(z_2,R_1+2t)\setminus B(z_1,R_1+2t)\big)\cup \big( B(z_1,R_1+2t)\setminus B(z_1,R_1)\big).
\]
Clearly, both set differences allow for upper bounds similar to those in \eqref{E:Inclusion1b} and \eqref{E:Inclusion3}.
Hence, it remains to derive an upper bound for the measure of the difference $\sup_{y_J \in B^{\otimes r}(x_J;t)}\|z_2-z_1\|$ when integrated w.r.t.\ $x_{J}$ where $J=J_1\cap J_2$ (to cover both situations). To this end, we study the difference of the centers of the smallest bounding ball when the index sets differ by exactly one element given $x_J$, $y_J$. Assume that $J=\{1,\ldots,k\}$, where $k\le d$ because the simplices of $x_{J\cup\{j_0\}}$, $y_{J\cup\{j_0\}}$ appear in the \v Cech filtration by assumption. Then we use a telescoping sum to obtain
\[
	\|z_2-z_1\| \le \sum_{\ell=1}^k \| z(x_i,y_j: i\le \ell < j) - z(x_i,y_j: i< \ell \le j) \|.
\]
Consider the summands on the right-hand side. When taking the supremum over $y_J\in B^{\otimes r}(x_J;t)$ for each of it, it is clearly sufficient to study the dominating term with $\ell = 1$.
We denote $d(x_1,U^*)$ the distance of $x_1$ to the affine linear space $U^*$ spanned by the elements in $\wt x = (x_2,\ldots,x_k)$. Put $\wt y = (y_2,\ldots,y_k)$ and consider
\begin{align}\label{E:IntegralSup}
\int_M   \sup_{y_1 \in B(x_1;t) } \sup_{\wt y \in B^{\otimes r}(\wt x;t) } \| z(x_1,\wt y) - z(y_1,\wt y) \| \ \diff x_1.
\end{align}
We apply Proposition~\ref{P:CechDiffCenters} with $\delta=1/2$: On the one hand provided $d(x_1,U^*) \ge (6\diam M +4) \cdot t$, we have
\[
	\sup_{y_1\in B(x_1;t)} \sup_{ \wt y \in B^{\otimes r}(\wt x;t)} \| z(x_1, \wt y ) - z(y_1, \wt y) \|	\le \frac{16\diam M + 2}{ d(x_1,U^*)} t.
	\]
On the other hand, both centers lie in $M$ by convexity, so can rely on the trivial upper bound $\| z(x_1,\wt y) - z(y_1,\wt y) \| \le \diam M$ in the case that $d(x_1,U^*) < (6\diam M +4) \cdot t$.

In order to bound \eqref{E:IntegralSup} from above, we set $c_1 \coloneqq 6\diam M +4$ and $c_2 \coloneqq 16\diam M + 2$ and study
\begin{align*}
	&\int_M \diam M \1{ d(x_1,U^*) < c_1 t }  \diff x_1 \\
	&+ \quad \int_M  \frac{ c_2 t}{d(x_1,U^*) } \1{ d(x_1,U^*) \ge c_1 t } \diff x_1 .
\end{align*}
Consider the first integral. As the affine linear space $U^*$ is defined by $k-1$ elements with $k\le d$, its dimension is at most $d-2$. Thus, by Lemma~\ref{L:AnnulusIntersectedManifold}
\begin{align*}
	\int_M  \diam M \1{ d(x_1,U^*) < c_1 t } \diff x_1 \lesssim (\diam M)^{1+(d-2)}  (c_1 t)^2 . 
\end{align*}
Next, we study the second integral. It is evident that distances $d(x_1,U^*) > 1$ make a contribution to the overall integral of order $t$. We introduce the additional condition $d(x_1,U^*) \in ((n+1)^{-1},n^{-1}]$ for some $n\in\N$. Then we obtain the upper bound
\begin{align*}
	& \sum_{n\in\N}\int_M   \frac{ c_2 t}{d(x_1,U^*) }  \1{ (c_1 t) \vee \frac{1}{n+1} < d(x_1,U^*) \le \frac{1}{n} } \diff x_1 \\
	&< \sum_{n\in\N} \int_M  \frac{ c_2 t}{(n+1)^{-1} }  \1{ \frac{1}{n+1} < d(x_1,U^*)\le \frac{1}{n} } \diff x_1 \\
	&=   \frac{ c_2 t}{(n+1)^{-1} }  \,\sum_{n\in\N } \, \lambda^d\Big(\Big\{ x_1\in M : \frac{1}{n+1} < d(x_1,U^*)\le \frac{1}{n} \Big \} \Big).
\end{align*}
Applying once more Lemma~\ref{L:AnnulusIntersectedManifold}, we see that this last series is of order
\[
	\sum_{n\in\N} \frac{t}{(n+1)^{-1}} \Big( \frac{1}{n^2} - \frac{1}{(n+1)^2} \Big) = \sum_{n\in\N} \frac{t (n+1)(2n+1)}{n^2(n+1)^2} \le \sum_{n\in\N} \frac{2t}{n^2} = \frac{\pi^2 t}{3}.
\]
Hence, we can conclude as in Case I the rate $t$, i.e., the exponent 1 (after integrating w.r.t.\ the remaining variables $\wt x$).

Subcase III: $J_1 \setminus J_2 \neq \emptyset \neq J_2\setminus J_1$. The setting in \eqref{E:Int2} implies
\begin{align}\label{E:IntCech1}
	\phi[J_2](x) \in ( \phi[J_1](x), \phi[J_1](x) + 2t]
\end{align}
and the setting in \eqref{E:Int3} simply implies
\begin{align}\label{E:IntCech2}
	\phi[J_2](x) = \phi[J_1](x).
\end{align}
Observe that we can assume in both cases that $J_2$ is minimal w.r.t.\ inclusion, viz., we have $\phi[J'](x)< \phi[J_2](x)$ whenever $J' \subset J_2$ is a strict subset and all $x_j$ lie on the circumsphere.
We write $R(x_{J})=\phi[J](x)$ for the filtration time of the tuple $x_J$. Assume that $j_0\in J_2 \setminus J_1$ and set $\wt J=J_2\setminus\{j_0\}$. We set
\[
	A_{\wt J} = \{x_{\wt J} : R(x_{\wt J}) > R(x_{J'}) \text{ for all strict subsets } J' \subset \wt J\}
\]
and note that due to the minimality of $J_2$, $x_{\wt J}\in A_{\wt J}$ whenever $x$ satisfies \eqref{E:IntCech1} or \eqref{E:IntCech2}.

We begin with \eqref{E:IntCech1} and set $a=R(x_{J_1})$, $b=R(x_{J_1})+2t$. Then we consider
\begin{align}\begin{split}\label{E:IntCech3}
		&\lambda( x_{j_0}\in M: R(x_{J_2}) \in (a,b) ) \cdot \1{x_{\wt J} \in A_{\wt J} }\\
		&\cdot \Big( \1{ R(x_{\wt J}) \le a }  +  \1{ a< R(x_{\wt J}) < b } \Big).
		\end{split}
\end{align}
Note that $R(x_{J_2})\in (a,b)$, $a<R(x_{\wt J})$ and $x_{\wt J} \in A_{\wt J}$ imply $R(x_{J_2})\in (R(x_{\wt J}),b)$.

Next, we rely on Proposition~\ref{P:ConditionalDensity} in order to determine this measure in \eqref{E:IntCech3} and integrate w.r.t. the remaining indices $\wt J$. We arrive at
\begin{align*}
	&\int_{A_{\wt J}}  \Big( \sqrt{ b^2 - R(x_{\wt J}) ^2} - \sqrt{ a^2 - R(x_{\wt J})^2}  \Big)  \1{ R(x_{\wt J}) \le a } \, \diff x_{\wt J}  \\
	&\qquad\qquad + \int_{A_{\wt J}}  \sqrt{ b^2 - R(x_{\wt J}) ^2} \1{ a < R(x_{\wt J}) <  b } \, \diff x_{\wt J} .
\end{align*}
In order to evaluate these two integrals, we rely on Divol and Polonik \cite{divol2019choice} who show that $R(x_{\wt J})$ when considered as a random variable on $\R_+$ admits a bounded density w.r.t.\ the Lebesgue measure. Thus, in order to derive an upper bound, we can calculate instead
\begin{align}\label{E:IntCech4}
	&\int_0^{a} \Big\{ \sqrt{ b^2 - u^2} - \sqrt{ a^2 - u^2} \Big\}  \, \diff u + \int_a^{b}  \sqrt{ b^2 - u^2} \, \diff u .
\end{align}
Remember the indefinite integral $\int \sqrt{V^2 -u^2} \, \diff u$ is $u \sqrt{V^2-u^2} + V^2 \arcsin(V^{-1} u)$. Thus, this last integral equals $\pi/4 \cdot (b^2-a^2) \le \pi \diam M/4 \cdot  t$. All in all, \eqref{E:IntCech1} is of order $t$.

We continue with \eqref{E:IntCech2}. We obtain
\begin{align}
	&\lambda( \{ x_{j_0} : R(x_{J_2}) = R(x_{J_1}) \} ) \1{x_{\wt J} \in A_{\wt J} }  \nonumber \\
	&=\lim_{\epsilon\downarrow 0} \lambda( \{ x_{j_0} : R(x_{J_2}) \in ( a - \epsilon, a + \epsilon) \1{x_{\wt J} \in A_{\wt J} },\label{E:IntCech5}
\end{align}
where $a=R(x_{J_1})$. Again, we can proceed as in \eqref{E:IntCech3} and study the subcases where either $R(x_{\wt J}) \le a-\epsilon$ or $a-\epsilon < R(x_{\wt J}) < a+\epsilon$ to arrive at \eqref{E:IntCech4} with $a$ being replaced by $a-\epsilon$ and $b$ by $a+\epsilon$.
Hence, taking the limit in \eqref{E:IntCech5} shows that the event underlying \eqref{E:IntCech2} is actually a null set.
\end{proof}

\begin{proposition}\label{P:CechDiffCenters}
Let $t,\delta\in (0,1)$ and $J\subset [r]$. Conditional on $x_J$, let $U^*$ be the smallest affine linear space containing $x_J$. Provided that
\[	
d(x_0,U^*) \ge (2\diam M +2)/(1-\delta) \cdot t,
\]
we have
\[
	\sup_{y_0\in B(x_0;t)} \sup_{ y_J\in B^{\otimes r}(x_J;t)} \| z(x_0, y_J) - z(y_0, y_J) \|	\le \frac{8\diam M + 1}{\delta d(x_0,U^*)} t.
	\]
\end{proposition}
\begin{proof} We split the proof in six steps. Let $y_J\in B^{\otimes r}(x_J;t)$.

Step 1. First consider the subset $\wt J$ of $J$ that is minimal w.r.t.\ inclusion such that $\phi[\wt J](y)=\phi[J](y)$. Then we have for the centers
\[	
	z(x_0,y_J) = z(x_0,y_{\wt J}) \text{ and } z(y_0,y_J) = z(y_0,y_{\wt J}).
\]
Denote $i\in \wt J$ an arbitrary but fixed index and consider the affine linear spaces
$$
	U\coloneqq U_{\wt J} \coloneqq y_i+span(y_j-y_i:j\in\wt J) \text{ and } U^*\coloneqq U^*_J\coloneqq x_i+span(x_j-x_i: j\in  J).
$$
We denote the orthogonal projection of some $w\in \R^d$ onto $U$ by $\Pi_{U}.w$.
Partition the domain $M$ conditionally on $y_{\wt J}$ in three sets 
\begin{align*}
	& B^*(y_{\wt J}),\, D(y_{\wt J}) \coloneqq \{ w \in M\setminus B^*(y_{\wt J}): \Pi_{U}.w \in B^*(y_{\wt J}) \},\\
	&\text{and } O(y_{\wt J}) \coloneqq M\setminus (B^*(y_{\wt J})\cup D(y_{\wt J})).
\end{align*}

Step 2. We partition the line segment $\overline{x_0 y_0}$ into the three subsegments (on the same line) $\overline{ x_0 x_0^r}$, $\overline{x_0^r y_0^r}$, $\overline{y_0^r y_0}$ such that $x_0^r$, $y_0^r$ lie both in $D(y_{\wt J})$ and the distances $\overline{x_0 x_0^r}$, $\overline{y_0 y_0^r}$ are minimal. Thus, $x_0$ (resp.\ $y_0$) can coincide with $x_0^r$ (resp.\ $y_0^r$) or $x_0^r$ can coincide with $y_0^r$ depending on the exact positions of $x_0$ and $y_0$.

Step 3. We study the difference of the centers of the smallest bounding balls belonging to $x_0^r$ and $y_0^r$. We abbreviate the latter in this step simply by $x_0$ resp.\ $y_0$. W.l.o.g.\ let $B^*(y_{\wt J})$ be centered at 0 with radius $R_0 = R_0(y_{\wt J})=\phi[\wt J](y)$. Then $\Pi_U.w$ is the orthogonal projection of $w$ onto the linear space $U$ and $w^{\perp_U}$ the projection onto the orthogonal complement. 

Set $\nu(x_0) = x_0^{\perp_U} / \| x_0^{\perp_U}\|$. We consider the center $z(x_0,y_{\wt J})$ of the smallest bounding sphere of $\{x_0,y_j:j\in\wt J\}$. The coordinate perpendicular to $U$ points in direction $\nu(x_0)$ and -- relying on the Pythagorean theorem -- has the value
\[
	h(x_0) = \frac{ \|x_0\|^2 - R_0^2 }{2 \|x_0^{\perp_U}\| }.
\]
Then $z(x_0,y_{\wt J}) = h(x_0) \cdot \nu(x_0)$. Similarly, $z(y_0,y_{\wt J}) = h(y_0)\cdot \nu(y_0)$. It remains to bound the difference $z(x_0,y_{\wt J}) - z(y_0,y_{\wt J})$ under the restriction that $\|x_0-y_0\| \le t$. On the one hand, we have 
\begin{align*}
	| \nu(x_0) - \nu(y_0) | \le \frac{ \big| \|x_0^{\perp_U}\| - \|y_0^{\perp_U}\| \big| }{ \|x_0^{\perp_U}\| } + \Bigg\| \frac{y_0^{\perp_U}}{\|x_0^{\perp_U}\|} - \frac{ y_0^{\perp_U}} {\|y_0^{\perp_U} \|} \Bigg\| \le 2 \frac{t}{\|x_0^{\perp_U}\|}.
\end{align*}
On the other hand,
\begin{align*}
	| h(x_0) - h(y_0) | &\le \frac{ \big| \|x_0^{\perp_U}\|^2 - \|y_0^{\perp_U}\|^2  \big|  } {2 \|x_0^{\perp_U}\| } + \Bigg| \frac{\|y_0\|^2 - R_0^2}{2 \|x_0^{\perp_U}\| } -  \frac{\|y_0\|^2 - R_0^2}{2 \|y_0^{\perp_U}\| }  \Bigg| \\
	&\le \frac{2 \|x_0\|+t}{2 \|x_0^{\perp_U}\|} \cdot t + \frac{\|y_0\|^2 - R_0^2}{2 \|x_0^{\perp_U}\|} \cdot \frac{2 | \|x_0^{\perp_U}\| - \|y_0^{\perp_U}\| | }{2\| x_0^{\perp_U}\| }  \\
	& \le \Big( \|x_0\| + \frac{t}{2} + |h(x_0)| \Big) \cdot \frac{t}{\|x_0^{\perp_U}\|} \le \Big(2\diam M + \frac{t}{2}\Big)\cdot \frac{t}{\|x_0^{\perp_U}\|}.
\end{align*}
Note that we have relied on the centering at the origin and on $|h(x_0)|\le \|x_0\|\le \diam M$ to derive the last inequality. Overall, we see that 
\[
	\| z(x_0,y_{\wt J}) - z(y_0, y_{\wt J}) \| \le \Big(4\diam M + \frac{t}{2}\Big) \frac{t}{ \| x_0^{\perp_U}\|}.
\]
It remains to bound $ \| x_0^{\perp_U}\| = d(x_0,U)$ from below by an expression involving $d(x_0,U^*)$. Consider the minimizer of $U\ni w\mapsto d(x_0,w)$ which we denote by $w^* = y_i + \sum_{j \neq i, j\in \wt J}\mu_j(y_j-y_i)$ for certain $\mu_j\in\R$. Then $d(w^*, M) \le \diam M$. Indeed, $x_0, w^*$ and any other point in $p\in U\cap M$ form a triangle such that the angle between $x_0-w^*$ and $p-w^*$ is $\pi/2$ and $\|x_0 - p\|\le \diam M$. Consequently,
\[
	d(w^*,M)^2 \le d(w^*,x_0)^2 \le \|x_0-w^*\|^2 + \|p-w^*\|^2 = \|x_0-p\|^2 \le (\diam M)^2.
\]
Then as $y_j = t\nu_j + x_j$ for some unit vector $\nu_j$ and as $ d(w^*,M) \ge \|w^*-x_0\|$, we have
\begin{align*}
	2\diam M \ge d(w^*,M) + \|x_0-y_i\| \ge \|w^* - y_i\| \ge  \Big\|\sum_{j \neq i, j\in \wt J}\mu_j(y_j-y_i) \Big\|,
\end{align*}
i.e., $ \| \sum_{j \neq i, j\in \wt J}\mu_j(y_j-y_i) \| \le 2\diam M$.
And it follows,
\begin{align*}
	\|x_0 - w^*\| &\ge \Big\|x_0 - \Big\{x_1 +  \sum_{j \neq i, j\in \wt J}\mu_j(x_j-x_i) \Big\}\Big \| - t \Big(1+\Big \|  \sum_{j \neq i, j\in \wt J}\mu_j(\nu_j-\nu_i) \Big\| \Big) \\
	&\ge \Big\|x_0 - \Big\{x_1 +  \sum_{j \neq i, j\in \wt J}\mu_j(x_j-x_i) \Big\}\Big \| - (2\diam M+1) t \\
	&\ge \min_{w\in U^*} d(x_0,w) - (2\diam M+1) t \\
	&= d(x_0,U^*) - (2\diam M+1)t.
\end{align*}
Putting the estimates together, we arrive at
\[
	\| z(x_0,y_{\wt J})-z(y_0, y_{\wt J}) \| \le \frac{(4\diam M + t/2) t}{d(x_0,U^*) - (2\diam M+1)t} \le \frac{(4\diam M + t/2) t}{\delta d(x_0,U^*)}
\]
because $d(x_0,U^*) \ge (2\diam M +2)/(1-\delta) \cdot t \ge (2\diam M +1)/(1-\delta) \cdot t$

Step 4. If $x_0\in B^*(y_{\wt J})$, then $x_0^r$ lies at the boundary of $B^*(y_{\wt J})$. Then plainly, $ z(x_0^r,y_{\wt J}) = 0 = z(x_0,y_{\wt J})$ Similar considerations apply if $y_0\in B(y_{\wt J})$.

Step 5. If $x_0\in O(y_{\wt J})$, then $x_0^r$ is in $\partial O(y_{\wt J}) \cap \partial D(y_{\wt J})$ and the smallest bounding sphere of $\{x_0, y_j: j\in \wt J\}$ is the same as that of $\{x_0,y_i: i\in I\}$ for a strict subset $I \subset \wt J$. 
Indeed the projection $\Pi_U.x_0$ is not contained in $B^*(y_{\wt J})$ (as otherwise $x_0\in B^*(y_{\wt J}) \cup D(y_{\wt J})$). Thus, the filtration time of $\{ \Pi_U.x_0, y_j: j\in\wt J\}$ is strictly larger than that of $y_{\wt J}$, however, all elements lie in the same subspace of minimal dimension. Hence, we can dispose of certain $y_j$.

Now, observe that the transition from $\{x_0,y_i: i\in I\}$ to $\{x^r_0,y_i: i\in I\}$ is just the setting in Step 3. If $U'$ is the affine space spanned by the $y_i$, $i\in I$, then $d(x_0,U') \ge d(x_0, U)$ because $U'$ is contained in $U$.

Similarly, if $y_0\in O(y_{\wt J})$, we find a strict subset $I$ of $\wt J$ achieving the same filtration time for $y_{\wt J}$ and $y_I$. Again, the transition from $\{y_0,y_i: i\in I\}$ to $\{y^r_0,y_i: i\in I\}$ is as in Step 3. However, in this case, the lower bound deteriorates by $-t$, we have
\[
	d(y_0,U') \ge d(y_0, U) \ge d(x_0,U)-t \ge d(x_0,U^*) - (2\diam M + 2)t.
\]

Step 6. In total, when performing the splits of the line segment $\overline{x_0 y_0}$ into $\overline{x_0 x_0^r}$, $\overline{x_0^r y_0^r}$, $\overline{y_0^r,y_0}$, we always have $0<t_1,t_2$ adding up to $t$ where $t_1$ contributes from Step 3 and $t_2$ contributes from Step 5. Note that we cannot have $x_0$ and $y_0$ both in $O(y_{\wt J})$ because initially we assumed that $U^*$ is minimal. So, in total we arrive at
\begin{align*}
	&\|z(x_0,y_J)-z(y_0,y_J)\|=\| z(x_0,y_{\wt J}) - z(y_0, y_{\wt J}) \| \le 2 \frac{(4\diam M + t/2)}{\delta d(x_0,U^*)} t.
\end{align*}
This completes the proof.
\end{proof}

\begin{lemma}\label{L:AnnulusIntersectedManifold}
Let $L$ be an affine subspace of $\R^d$ of dimension $k\le d-2$ such that $L\cap M\neq \emptyset$ and $0\le a < b$. Then
\begin{align}\begin{split}\label{E:AnnulusIntersectedManifold1}
	&\lambda^d\Big(\Big\{ z \in M : a < d(z,L)\le b \Big \} \Big) \\
	&\le \frac{\pi^{(d-k)/2}}{\Gamma((d-k)/2+1)} (b^{d-k} - a^{d-k}) \cdot \frac{\pi^{k/2}}{\Gamma(k/2+1)} (\sqrt{d}\diam M)^k.
\end{split}
\end{align}
\end{lemma}
\begin{proof}
Observe that upon replacing $M$ by a cube of side length $\diam M$ and replacing $L$ by the linear space $L^* = \{ z \in \R^d : z_{k+1} = \ldots =z_d = 0\}$, we obtain
\[
	\lambda^d(\{ z \in M : a < d(z,L)\le b  \} ) \le \lambda^d (\{ z \in [0,\diam M]^d : a < d(z,L^*)\le b \} ). 
\]
Clearly, as in this case $d(z,L^*)^2 = z_{k+1}^2+\ldots+z_d^2$, the right-hand side is dominated by
\begin{align*}
	&\lambda^d(\{z\in \R^d: \|z\|\le \sqrt{d}\diam M \text{ and } a^2 \le z_{k+1}^2+\ldots+z_d^2 \le b^2 \} ) \\
	&\le \lambda^{k}\Big(\Big\{(z_i: i\le k): \sum_{i=1}^k z_i^2 \le d(\diam M)^2 \Big\}\Big) \\
	&\quad \qquad \cdot \, \lambda^{d-k}(\{ (z_i: i > k) : a^2 \le z_{k+1}^2+\ldots+z_d^2 \le b^2	\})
\end{align*}
which in turn is at most the product given in \eqref{E:AnnulusIntersectedManifold1}.
\end{proof}

\begin{proposition}[Conditional density of the \v Cech filtration function]\label{P:ConditionalDensity}
Let $\wt J = \{1,\ldots,k\}$ and $J=\{j_0\}\cup \wt J$ with $k<j_0\le d$.
Let $x_{\wt J }\in M^k$ such that $\phi[J'](x)< \phi[\wt J](x)=:R(\wt x)$ for all strict subsets of $\wt J$. Let $a,b\in\R_+$ such that $R(x_{\wt J}) \le a < b$ and consider $R(x_J) \coloneqq \phi[J](x)$. We have
\[
	\lambda^{d}\Big(\{ x_{j_0} \in M: R(x_J) \in (a,b) \} ) \le C \Big( \sqrt{ b^2 - R(x_{\wt J})^2 } - \sqrt{ a^2 - R(x_{\wt J})^2 } \Big)
\]
for a universal $C\in\R_+$ which is independent of $\wt J$ and $x_{\wt J}$.
\end{proposition}
\begin{remark}
It is important to note that the estimate is only valid under conditions that guarantee $\phi[J](x) > \phi[\wt J](x)$ because the measure of $\{x_{j_0} \in M: \phi[J](x) = \phi[\wt J](x)\}$ is positive given general $x_{\wt J}$.
\end{remark}
\begin{proof}[Proof of Proposition~\ref{P:ConditionalDensity}]
We build on existing results of Divol and Polonik \cite{divol2019choice} that lay out in detail the geometric situation in this case, however, we adopt the techniques to our context. We define $U \coloneqq span(x_{j} - x_1 : j\in \wt J)$. Observe that all $x_j$ which span $U$ lie on the circumsphere of $x_{\wt J}$. We assume that the circumcenter $\wt z \coloneqq z(x_{\wt J})$ of $x_{\wt J}$ is at the origin; if the circumcenter does not coincide with the origin, we have to subtract it from all $x_j$. Thus, $U$ is a proper linear space and has an orthogonal complement $U^\perp$. Then, define
$$
	 V(\theta) \coloneqq U + \theta \R \text{ for } \theta\in U^{\perp} \text{ with } \|\theta\|=1. 
$$
Next, we choose an arbitrary but fixed $\theta_0\in U^{\perp}$ with norm 1, define $V_0 = V(\theta_0)$. Here the choice of $\theta_0$ solely depends on $x_{\wt J}$. Moreover, define the unit sphere in $U^\perp$ (resp.\ in $V_0$) by $S_{U^\perp}$ (resp. $S_{V_0}$).
Define for $\theta\in U^\perp$ the isometry from $V_0$ to $V(\theta)$ by
\[
	T(\theta).v \coloneqq \Pi_U.v + \langle v,\theta_0\rangle \theta,
\]
where $\Pi_U$ is the orthogonal projection onto $U$. Moreover, we decompose a vector $v\in\R^d$ orthogonally into $\Pi_U.v$ and $v^{\perp_U}$.
Define conditionally on $\wt R \coloneqq R(x_{\wt J})$
\[
	\Phi\colon \R_+\times S_{U^\perp}\times S_{V_0} \to \R^d, (t,\theta,\nu)\mapsto t\theta + \sqrt{t^2+\wt R^2} \cdot T(\theta).\nu
\]
and put for $a<b$
\[
	M_{a,b} \coloneqq \{ s\ge 0 : a^2 \le s^2 + \wt R ^2 \le b^2 \}.
\]
More details on the construction are given in Divol and Polonik~\cite{divol2019choice}. The important point is at this stage that due to the minimality of $J$ the filtration time $R\coloneqq R(x_{J})$ satisfies the Pythogorean theorem: If $z$ is the circumcenter of $x_{J}$ and $\wt z=0$ that of $x_{\wt J}$, then $z-\wt z$ is orthogonal to $U$ and points in the direction of $(x_{j_0}-\wt z)^{\perp_U}$, i.e., the angle between $z-\wt z$ and $(x_{j_0}-\wt z)^{\perp_U}$ is zero. Thus, 
\[
		\| (x_{j_0} - \wt z) - (z-\wt z)\|^2 = \|x_{j_0} - z \|^2 = R^2 = \|z - \wt z\|^2 + \wt R^2.
\]
Hence, $R\in [a,b]$ is true if and only if $a^2 \le t^2 + \wt R^2 \le b^2$ for $t \coloneqq\|z - \wt z\|$.

Set $t \coloneqq \|z-\wt z\|$ and $\theta^* \coloneqq (z-\wt z)/\|z-\wt z\|$. Then the unit vector $\theta^*$ is in the orthogonal complement of $U=span(x_{\wt J})$ as required. Moreover, choose $\nu \in S_{V_0}$ of unit norm such that
\[
	(x_{j_0}-z) / R  = \Pi_U.\nu + \langle \nu, \theta_0 \rangle \theta^*.
\]
Indeed such a $\nu$ exists: Observe that $(x_{j_0}-z) / R$ has unit norm. Set $s \coloneqq \sqrt{ 1 - \| \Pi_U.(x_{j_0}-z) / R \|^2 }$, note that the argument inside the square root is non-negative. Define
\[
		\nu \coloneqq \Pi_U.(x_{j_0}-z) / R + \theta_0 s
\]
Then $\langle \nu, \theta_0\rangle = s$. So,
\begin{align*}
	T(\theta^*).\nu &= \Pi_U.\nu + \langle \nu, \theta_0 \rangle \theta^* = \Pi_U.(x_{j_0}-z) / R + s \theta^* \\
	&= \frac{ \Pi_U.(x_{j_0}-z) }{ R } + \sqrt{ 1 - \| \Pi_U.(x_{j_0}-z) / R \|^2 } \cdot \frac{z-\wt z}{\|z-\wt z\|\ }\\
	&= \frac{ x_{j_0} - z }{R} = \frac{ (x_{j_0} - \wt z) - (z-\wt z) }{R}.
\end{align*}
Or equivalently, $x_{j_0}$ satisfies in this case (with the quantities as just defined)
\[
		x_{j_0} - \wt z = t\theta^* + \sqrt{t^2+\wt R^2} \cdot T(\theta^*).\nu.
\]
In particular, $x_{j_0} - \wt z \in \Phi( M_{a,b} \times S_{U^\perp} \times S_{V_0})$.

With these preparations we come back to \eqref{E:IntCech1} and \eqref{E:IntCech2}. Then given the choice of $\theta_0$, which solely depends on $x_{\wt J}$, we arrive at
\begin{align*}
	&R(x_{J}) \in (a , b) \text{ and } a \ge R(x_{\wt J}) \text{ implies } x_{j_0} \in \Phi( M_{a,b} \times S_{span(x_{\wt J})^\perp} \times S_{V(\theta_0)}) .
\end{align*}
Divol and Polonik show that under these conditions and conditional on $x_{\wt J}$ (rephrased in our setting) 
\begin{align*}
	& \lambda^d\Big( x_{j_0} \in M: x_{j_0} \in \Phi( M_{a,b} \times S_{span(x_{\wt J})^\perp} \times S_{V_0} ) \Big) \\
	&\le C  \int_{a}^{b} \frac{s^d}{\sqrt{s^2-R(x_{\wt J})^2 }} \diff s
\end{align*}
for a suitable constant $C\in\R_+$ detailed in this paper. Next, we use that $a,b \le \frac{1}{2}\diam M < \infty$ to find that this last integral is of order
\begin{align*}
	\lesssim \left(\frac{\diam M}{2} \right)^{d-1} \int_{a}^{b} \frac{s}{\sqrt{s^2-R(x_{\wt J})^2 }} \diff s &\lesssim \sqrt{ b^2 - R(x_{\wt J})^2 } - \sqrt{ a^2 - R(x_{\wt J})^2 } 
\end{align*}
because $d\ge 2$.
\end{proof}

\begin{proof}[Proof of Theorem~\ref{T:GradientCondition}]
We split the proof in two main steps.

{\it Step 1.} We start with the situation in \eqref{E:Int2}. Then relying on the mean value theorem and the convexity of $M$ there are $\theta_1,\theta_2\in [0,1]$ which satisfy
\begin{align}\begin{split}\label{E:Int1}
	&\phi[J_1](x) + \langle \nabla \phi[J_1](x+\theta_1(y-x)), y-x \rangle = \phi[J_1](y) \\
	&\ge \phi[J_2](y) = \phi[J_2](x) + \langle \nabla \phi[J_2](x+\theta_2(y-x)), y-x \rangle .
\end{split}\end{align}
Thus, relying on the Cauchy-Schwarz-inequality and $\|x-y\|^2 \le r t^2$, we obtain
\begin{align}\label{E:Int4}
	&\phi[J_1](x) < \phi[J_2](x) \le \phi[J_1](x) + 2 c^* \sqrt{r} t.
\end{align}
Note that if the \v Cech filtration is applied, we can utilize \eqref{E:GradientCech1} in Remark~\ref{R:GradientCech} (instead of the mean value theorem) to tighten the interval as follows
\begin{align*}
	\phi[J_2](x) \in (\phi[J_1](x), \phi[J_1](x) + 2 t ].
\end{align*}
If the Vietoris-Rips filtration is applied, we utilize Remark~\ref{R:GradientVR}: $\nabla\phi[J](x)$ has norm $\sqrt{2}$ and can be identified with the partial gradients $\nabla^1\phi[J](x) ,\ldots,\nabla^r\phi[J](x)$ two of which are nonzero. Hence, relying on \eqref{E:Int1}, we obtain
\begin{align*}
	\phi[J_2](x) \in (\phi[J_1](x), \phi[J_1](x) + 2 \sqrt{2} t ].
\end{align*}

{\it Step 2.} It remains to consider the situation in \eqref{E:Int3} where we can assume w.l.o.g. that $J_1\cap J_2\neq\emptyset$, see Remark~\ref{R:Rho}.
Observe that $J_1\subseteq J_2$ is not possible because $\phi[J_1](y)> \phi[J_2](y)$, Moreover, remember $J_2$ is minimal with respect to inclusion but we can have $J_2\subset J_1$.

Firstly, if both $J_1\setminus J_2 \neq\emptyset  \neq J_2\setminus J_1$, then relying once more on the mean-value theorem the situation in \eqref{E:Int3} is included in 
\begin{align}\label{E:Int5}
	\phi[J_2](x) \in [\phi[J_1](x), \phi[J_1](x) + 2c^*\sqrt{r}t).
\end{align}
Plainly, if the \v Cech filtration function is applied the right endpoint of the interval can be replaced by $2t$; similarly we can replace it by $2\sqrt{2}t$ for the Vietoris-Rips filtration function.

Hence, \eqref{E:Int4} and \eqref{E:Int5} show \eqref{E:GradientCondition1} for any two different sets $J_1,J_2$ that are not subsets of each other.

Secondly, if $J_2\subset J_1$ and if \eqref{E:Int3} is satisfied, we conclude once more (similar to \eqref{E:Int1}) the inequalities
\[
	\phi[J_2](y) < \phi[J_1](y) \le \phi[J_1](x) +c^*\sqrt{r}t = \phi[J_2](x) +c^*\sqrt{r}t  \le \phi[J_2](y) +2c^*\sqrt{r}t.
\]
Hence, $\phi[J_1](y) \in (\phi[J_2](y) , \phi[J_2](y) +2c^*\sqrt{r}t]$.

Since the configuration $x$ satisfies $d(x,y)\le t$ for some $y$, we arrive at
\[
	x\in 	\{ y\in M^r: \phi[J_1](y) \in (\phi[J_2](y) , \phi[J_2](y) +2c^*\sqrt{r}t] \}^{(t)}.
\] 
This shows \eqref{E:GradientCondition2} in the case of an inclusion between $J_1$ and $J_2$.

All in all, this yields the assertion because the number of admissible sets $J_1,J_2$ is bounded above by $\sum_{k_1=1}^r \sum_{k_2=1}^r \sum_{j=0}^{k_1\wedge k_2} \binom{r}{k_1} \binom{r-k_1}{k_2-j}$.
\end{proof}

\subsection{$L^p$-$m$-approximable features}
\begin{proof}[Proof of Theorem~\ref{T:LpmApprox}]
Let $x, y \subseteq M^r$ be two generic sets consisting of $r$ elements. We consider a map $\phi[J](x) \to \phi[J](y)$ between filtration times. 
Then given that $x$ and $y$ are sufficiently close in the sense that $y \in B^{\otimes r}(x;t)$ and $\rho(x) > t$ for some $t>0$, the following map is a 
bijection between persistence diagrams $\pd(x)$ and $\pd(y)$ (of any admissible dimension $k$)
$$
	\gamma_N \colon \pd(x) \to \pd(y),\quad (b,d) = (\phi[J_b](x), \phi[J_d](x) ) \mapsto (\phi[J_b](y), \phi[J_d](y) ),
$$
where by convention for each $(b,d)\in \pd(x)$ there are $J_b,J_d \subseteq [r]$ such that $(b,d) = (\phi[J_b](x), \phi[J_d](x) )$.

Indeed, $\gamma_N$ is a bijection because by the definition of $\rho$ there are certain sets $J_i \subseteq [r]$, $1\le i \le L$, such that the two filtrations are given by
\begin{align*}
		&\phi[J_1](x) < \phi[J_2](x) < \ldots < \phi[J_L](x), \\
		&\phi[J_1](y) < \phi[J_2](y) < \ldots < \phi[J_L](y)
\end{align*}
for all $y \in B^{\otimes r}(x;t)$ provided $\rho(x)>t$. Consequently, the persistence diagrams can be identified via $\gamma_N$ in this case.

Next, let $Z_{k,0}^{(m)}$ be the representation of the persistence diagram obtained in the same fashion as $Z_{k,0}$ in  \eqref{E:RepresentationPD1} when replacing the data $\cX_0$ by $\cmX_0 = (X_{0,1}^{(m)},\ldots,X_{0,r}^{(m)} )$. We abbreviate in the sequel $Z_{k,0}$ by $Z_0$ and $Z^{(m)}_{k,0}$ by $Z^{(m)}_0 $. The number of features in $Z_0$ is bounded above by $N_k$.

Define the event
$$
	A_m =  \{ \cmX_0 \in  B^{\otimes r}(\cX_0;t) \text{ for a } t\in\Q, \  t < \rho(\cX_0)\}.
$$
Then it follows from the above considerations that given $A_m$ there is a one-to-one correspondence of the features listed in $Z_{k,0}$ and $Z_{k,0}^{(m)}$, i.e., if $(b,d)$ corresponds to the entries $\ell$ and $\ell+1$ in $Z_{k,0}$ and $(b^{(m)},d^{(m)})$ to the entries $\ell$ and $\ell+1$ in $Z_{k,0}^{(m)}$, then there are sets $J,\wt J \subset [r]$ such that
\begin{align}\begin{split}\label{E:LpmApprox2}
	&(b,d) = (\phi[J](\cX_0), \phi[\wt J](\cX_0)) \\
	&\qquad \text{ and } (b^{(m)},d^{(m)}) = (\phi[J](\cmX_0), \phi[\wt J](\cmX_0)).
\end{split}
\end{align}
Hence, on $A_m$ the difference $Z_{k,0} - Z_{k,0}^{(m)}$ can be interpreted as differences of filtration times w.r.t.\ the same sets.

Before we come to the final step, note that we have for the filtration time of a simplex $J\subseteq [r]$ in $\cX_0$, resp.\ $\cmX_0$, the universal upper bound 
\begin{align*}
	\max\{ \phi[J](\cX_0), \phi[J](\cmX_0) \} \le T.
\end{align*}
Furthermore, we set again $c^* \coloneqq \max_{J\subseteq [r]}\operatorname{ess \ sup}_{x\in M^r} \|\nabla \phi[J](x)\| < \infty$, where $\nabla \phi[J]$ denotes the gradient of $\phi[J]$.

We come to the final step. Let $p\ge 1$, then plainly
\begin{align}\label{E:LrmApprox3}
	\begin{split}
	 \E[ \| Z_0 - Z^{(m)}_0 \|^p]^{1/p} 
	&\le \E[ \| Z_0 - Z^{(m)}_0 \|^p \ \mathds{1}_{A_m} ]^{1/p} \\
	& \quad + \E[ \| Z_0 - Z^{(m)}_0 \|^p \ \mathds{1}_{A_m^c }] ^{1/p}.
	\end{split}
\end{align}
Relying on the relations in \eqref{E:LpmApprox2}, we have on $A_m$ the $\omega$-wise equalities
\begin{align*}
	&\| Z_0 - Z^{(m)}_0 \|^2 = \sum_{j =1}^{N_k} (d_{0,j} - d_{0,j}^{(m)})^2 + ( b_{0,j} - b_{0,j}^{(m)} )^2 \\
	&=  \sum_{j =1}^{N_k} ( \phi[\wt J_j](\cX_0) - \phi[\wt J_j](\cmX_0))^2 + \sum_{j =1}^{N_k} ( \phi[J_j](\cX_0) - \phi[J_j](\cmX_0))^2
\end{align*}
for suitable sets $J_j, \wt J_j \subseteq [r]$.

Next, relying on the assumptions on the filtration function, we find with the mean value theorem that
$$
	| \phi[J](\cX_0) - \phi[J](\cmX_0) |^2 \le (c^*)^2 \sum_{i\in J} \| X_{0,i} - X_{0,i}^{(m)} \|^2.
$$
Therefore the first term on the right-hand side of \eqref{E:LrmApprox3} is bounded by
\begin{align*}
	\E[ \| Z_0 - Z^{(m)}_0 \|^p \ \mathds{1}_{A_m} ]^{1/p} &\le \E[ (2 N_k \cdot (c^*)^2 \sum_{i=1}^r  \| X_{0,i} - X_{0,i}^{(m)} \|^2 )^{p/2} ]^{1/p} \\
	&\le \sqrt{2  N_k} c^* \cdot  \E[ (\sum_{i=1}^{r} \| X_{0,i} - X_{0,i}^{(m)} \|)^p]^{1/p} \\
	&\le \sqrt{2  N_k} c^* \cdot \sum_{i=1}^{r} \E[  \| X_{0,i} - X_{0,i}^{(m)} \|^p]^{1/p} .
\end{align*}
In order to bound the second term on the right hand side of \eqref{E:LrmApprox3} note that 
\begin{align*}
	\E[ \mathds{1}_{A_m^c}] \le 
	\p\big( \| X_{0,i} - X_{0,i}^{(m)} \| >  \rho(\cX_0) ~\text{for some i} \big) \le \sum_{i=1}^{r} \p\big( \| X_{0,i} - X_{0,i}^{(m)} \| >  \rho(\cX_0) \big).
\end{align*}
Hence, using the crude estimate $\| Z_0 - Z^{(m)}_0 \| \leq 2 N_k T$,
we obtain 
\begin{align*}
	\E[ \| Z_0 - Z^{(m)}_0 \|^p \  \mathds{1}_{A_m^c} ]^{1/p}	\le 2 N_k T \cdot \left( \sum_{i=1}^{r} \p\big( \| X_{0,i} - X_{0,i}^{(m)} \| >  \rho(\cX_0) \big) \right)^{1/p}.
\end{align*}
Let $u_m>0$. The probabilities inside the last factor are at most
\begin{align*}
	\p( \| X_{0,i} - X_{0,i}^{(m)} \| >  \rho(\cX_0))
	&\le \p( \| X_{0,i} - X_{0,i}^{(m)} \| > u_m ) + \p( \rho(\cX_0) < u_m )\\
	&\le u_m^{-p} \ \E[\| X_{0,i} - X_{0,i}^{(m)} \|^p] + C u_m^{\alpha}
\end{align*}
for a suitable $\alpha>0$ and $C\in\R_+$ by the condition in \eqref{E:Rho-Assumption} because $\cX_0$ admits a bounded density on $M^r$.
Consequently, 
\begin{align*}
	&\left( \sum_{i=1}^{r} \p( \| X_{0,i} - X_{0,i}^{(m)} \| >  \rho(\cX_0)) \right)^{1/p} \lesssim \sum_{i=1}^{r} u_m^{-1} \ \E[ \| X_{0,i} - X_{0,i}^{(m)} \|^p ]^{1/p} 
	+ u_m^{\alpha/p}.
\end{align*}
Employing these bounds in \eqref{E:LrmApprox3} yields 
\begin{align*}
	& \E[ \| Z_0 - Z^{(m)}_0 \|^p]^{1/p}\\
	&\lesssim  \sum_{i=1}^{r}  \E[ \| X_{0,i} - X_{0,i}^{(m)} \|^p]^{1/p} +   \sum_{i=1}^{r} u_m^{-1}  \E[ \| X_{0,i} - X_{0,i}^{(m)} \|^p ]^{1/p} 
	+ u_m^{\alpha/p} \\
	&= (1+u_m^{-1}) \gamma_p(m) + u_m^{\alpha/p},
\end{align*}
where $\gamma_p(m) = \sum_{i=1}^r \E[ \| X_{0,i} - X_{0,i}^{(m)} \|^p]^{1/p} $. So, the choice $u_m = \gamma_{p}(m)^{p/(\alpha+p)}$ yields $ \E[ \| Z_0 - Z^{(m)}_0 \|^p]^{1/p} \lesssim \gamma_p(m)^{\alpha/(\alpha+p)}$ and the right hand side is summable with respect to $m$ by Assumption~\ref{Assumption:Lrm-Condition-PointProcess}.
This completes the proof.
\end{proof}

\begin{proof}[Proof of Theorem~\ref{Theorem:FCLT}]
By Theorem \ref{T:LpmApprox} the process $(Z_{k,j})_j$	is $L^2$-$m$-approximable and so is $(f(Z_{k,j}))_j$ by the Lipschitz-continuity of $f$. 
Then the assertion follows from Theorem A.1 in Aue et al.\ \cite{aue2009break}.
\end{proof}

\begin{proof}[Proof of Corollary~\ref{Cor:FCLT}]
This is a direct consequence of Theorem \ref{Theorem:FCLT} and the continuous mapping theorem. We have for $t\in [0,1]$ 
\begin{align*}
	&\frac{1}{\sqrt{n}} \left( \sum_{j=1}^{[nt]} f(Z_{k,j}) - t \sum_{j=1}^n f(Z_{k,j}) \right) \\
	=& \frac{1}{\sqrt{n}} \sum_{j=1}^{[nt]} (f(Z_{k,j}) - \E f(Z_{k,j})) - \frac{t}{\sqrt{n}} \sum_{j=1}^n (f(Z_{k,j}) - \E f(Z_{k,j})) \\
	&\quad + \frac{[nt] - nt}{\sqrt{n}} \E f(Z_{k,j}).
\end{align*}
When considered as functions on $[0,1]$, the first two terms converge weakly to $B_{\Gamma_k} = W_{\Gamma_k} - \cdot \times W_{\Gamma_k}(1) $ in $D_\ell([0,1])$. The last term is of order $n^{-1/2}$ uniformly in $t\in [0,1]$.
\end{proof}

\subsection{Asymptotics under $H_1$}
\begin{proof}[Proof of Theorem~\ref{T:AsympH1}]
It is known that any stationary $L^2$-$m$-approximable sequence $(Y_i)_{i\in\Z}$ satisfies the strong law of large numbers, i.e.,
$$
	n^{-1} \sum_{i=1}^n Y_i \rightarrow \E[Y_1] \quad a.s. \quad (n\to\infty),
$$
see Theorem~\ref{L:SLLN} or for more context Rademacher and Krebs \cite[Lemma 3.2]{rademacher2023twosample}.
Consequently, by Theorem~\ref{L:USLLN}, it follows that for each $0\le a\le b\le 1$ and $\epsilon>0$
\begin{align*}
	&\p\Big( \sup_{t\in [a,b] } \Big| \Big(n^{-1} \sum_{i=1}^{ \floor{nt}} Y_i \Big)- t \E[Y_1] \Big| > \epsilon \Big) \\
	&\le \p\Big( \max_{1\le \ell \le n}  \frac{\ell+1}{n} \ \Big| \ell^{-1} \sum_{i=1}^{ \ell} Y_i -  \E[Y_1] \Big| \ge \frac{\epsilon}{2} \Big) + \1{ \frac{ |\E[Y_1]|}{n} \ge \frac{\epsilon}{2} } = o(1).
\end{align*}
Now, this last result can be applied to the pre-change sequence $(f( Z_{k,j}) )_j$ on the interval $t\in [0,\theta]$ and the post-change sequence $(f( Z^*_{k,j}) )_j$ on the interval $t\in (\theta,1]$ seperately. We omit the details.
\end{proof}

\begin{proof}[Proof of Theorem~\ref{T:HatTheta}]
The main factor decomposes as follows
\begin{align}\label{E:HatTheta1}
	S_v^T \wt\Sigma_n S_v &= \frac{1}{n} \cE_v^T \wt\Sigma_n \cE_v + 2 \frac{b_v}{n} \cE_v^T \wt\Sigma_n \Delta + \frac{b_v^2}{n} \Delta^T \wt\Sigma_n \Delta \quad (1\le v\le n),
\end{align}
for $b_v \coloneqq v(1-\theta)\1{v\le v^*} + \theta(n-v)\1{v>v^*}$.
We inspect the three terms in \eqref{E:HatTheta1} separately.
Firstly, it follows from \eqref{E:SLLNInno} that
\begin{align}\label{E:HatTheta4}
	\max_{v\le n} \| \cE_v \| \le \max_{v\le n} ~ \Big\| \Big(1-\frac{v}{n} \Big) \sum_{j=1}^{v} e_j - \frac{v}{n} \sum_{j=v+1}^n e_j \Big\| = o_{a.s.}(a_n).
\end{align}
Indeed, let $\delta>0$ be arbitrary but fixed. Since $ \sum_{j=1}^{n} e_j(\omega) / a_n \to 0$ for all $\omega\in\Omega_0$ measurable with $\p(\Omega_0)=1$, we find for each $\omega \in \Omega_0$ an $N_0(\omega)\in\N$ and a $C(\omega)\in\R_+$ such that
\begin{align*}
	&a_v^{-1} \sum_{j=1}^v e_j(\omega) \le \delta \quad \forall v\ge N_0(\omega) \text{ and } \sup_{v\in\N} a_v^{-1} \sum_{j=1}^v e_j(\omega) \le C(\omega).
\end{align*}
Combining these two estimates shows $\limsup_{n\to\infty} \max_{v\le n} |\sum_{j=1}^v e_j(\omega)| / a_n \le \delta$ and since $\delta$ is arbitrary this yields \eqref{E:HatTheta4}.
In particular, we obtain for the first two terms on the right-hand side of \eqref{E:HatTheta1}
\begin{align*}
	& \max_{v\le n} ~ n^{-1} \cE_v^T \wt\Sigma_n \cE_v = o_{a.s.} ( \| \wt\Sigma_n \| n^{-1} a_n^2  ) \text{ and }\\
	& \max_{v\le n} \frac{b_v}{n} | \cE_v^T \wt\Sigma_n \Delta | = o_{a.s.}( \|\wt\Sigma_n\| a_n ).
\end{align*}
Moreover,
$$
	\max_{v\le n} \frac{ b_v^2}{n} \Delta^T \wt\Sigma_n \Delta = n \theta^2(1-\theta)^2 \Delta^T \wt\Sigma_n \Delta + O( n^{-2} \|\wt\Sigma_n \|).
$$
Hence, \eqref{E:HatTheta1} equals
\begin{align}
	S_v^T \wt\Sigma_n S_v &= \frac{b_v^2}{n} \Delta^T\wt\Sigma_n\Delta + R_{n,v} \label{E:HatTheta2} \\
	&= n \ S^*\Big( \frac{v}{n}\Big)^T \wt\Sigma_n S^*\Big( \frac{v}{n}\Big) + R_{n,v}, \nonumber
\end{align}
where $\max_{v\le n} |R_{n,v}| = o_{a.s.}( \|\wt\Sigma_n\| a_n +  \|\wt\Sigma_n\| n^{-1} a_n^2)$. Consider  \eqref{E:HatTheta2}; we show that maximizing the left-hand side ultimately corresponds to maximizing the right-hand side on the set $\Omega_0$. To this end, set $\gamma_n \coloneqq a_n/n = o(1)$. We distinguish two cases $v\le v^*$ and $v>v^*$.

If $v\le v^*$, we have
\begin{align}\label{E:HatTheta3}
	\Big| \frac{b_{v^*}^2}{n} \Delta^T\wt\Sigma_n\Delta - \frac{b_{v^* - \floor{n\gamma_n}}^2}{n} \Delta^T\wt\Sigma_n\Delta \Big| \asymp  n \big(\gamma_n + O(\gamma_n^2) \big) \Delta^T\wt\Sigma_n\Delta .
\end{align}
In particular, the remainder $R_{n,v}$ is negligible when compared to this last difference because 
$$
	\frac{ \max_{v\le n} |R_{n,v}| }{  n \gamma_n  \Delta^T\wt\Sigma_n\Delta  }  = o_{a.s.}\Big( \frac{ \|\wt\Sigma_n \| (a_n + n^{-1} a_n^2) }{ n \gamma_n \Delta^T\wt\Sigma_n\Delta  } \Big) = o_{a.s.}(1).
$$
Consider a maximizer of \eqref{E:HatTheta2} and denote it by $\wh m$. If $\wh m < v^* - \floor{n\gamma_n}$, then we have for all $v$ between $v^* - \floor{n\gamma_n}$ and $v^*$
$$
	\frac{b_{\wh m}^2}{n} \Delta^T\wt\Sigma_n\Delta + R_{n,\wh m} \ge \frac{b_v^2}{n} \Delta^T\wt\Sigma_n\Delta + R_{n,v}
$$
which implies
\begin{align*}
	 0 &\le \frac{b_v^2 - b_{\wh m}^2}{n}   \Delta^T\wt\Sigma_n\Delta \\
	& \le 2 \max_{v\le n} |R_{n,v}| = o_{a.s.}\big( \|\wt\Sigma_n\| (a_n + n^{-1} a_n^2 ) \big) \qquad (v^* - \floor{n\gamma_n} \le v \le v^* ).
\end{align*}
Consequently, since the map $v\mapsto \frac{b_v^2}{n} \Delta^T\wt\Sigma_n\Delta$ is increasing for $1\le v\le v^*$, we obtain a contradiction to \eqref{E:HatTheta3} because this last bound implies that
$$
	\frac{b_{v^*}^2 - b_{v^* - \floor{n\gamma_n}}^2}{n}   \Delta^T\wt\Sigma_n\Delta  =  o_{a.s.}\big( \|\wt\Sigma_n\| (a_n + n^{-1} a_n^2 ) \big),
$$
too. This shows that all $\wh m$ which maximize \eqref{E:HatTheta2} satisfy $\wh m \ge v^* - \floor{n\gamma_n}$ for $n$ sufficiently large with probability 1.

In a similar spirit, one verifies that these $\wh m$ satisfy $\wh m\le v^* + \floor{n\gamma_n}$ for $n$ sufficiently large with probability 1. So any maximizer ultimately lies in the interval $[v^*- \floor{n\gamma_n}, v^* + \floor{n\gamma_n}]$ with probability 1.
Thus, any sequence $(\wh\theta_n)_n$ satisfies ultimately
$$
	|\wh\theta_n - \theta| \le \Big|\frac{\wh m}{n} - \frac{v^*}{n}\Big| + \Big|\frac{\floor{\theta n}}{n} - \frac{\theta n}{n}\Big| = O_{a.s.}(\gamma_n)+ O(n^{-1}) = O_{a.s.}( n^{-1} a_n).
$$
This completes the proof.
\end{proof}

\appendix
\section{External results}

\begin{theorem}[Theorem A.1 in Aue et al. \cite{aue2009break}]
Let $d,d'>0$. Let $(\epsilon_j:j\in\Z)$ be a sequence of iid random variables with values in $\R^{d'}$ and let $f\colon\R^{d'\times\infty}\to\R^d$ be a measurable function.
Define the $d$-dimensional random process $(y_j: j\in \Z)$ by
$$
	y_j = f(\epsilon_j,\epsilon_{j-1},\ldots), \qquad j\in\Z.
$$
Assume that $\E[ y_j ] =0$ and $\E[\|y_j\|^2 ] <\infty$. Suppose further that for any $m\ge 1$, the vectors $y_0^{(m)}$ are defined via
$$
	y_j^{(m)} = f^{(m)} (\epsilon_j,\ldots,\epsilon_{j-m}),\qquad j\in\Z,
$$
with measurable functions $f^{(m)}\colon\R^{d'\times(m+1)}\to\R^d$.
Assume that
$$
	\sum_{m\ge 1} \E[ \| y_0 - y_0^{(m)} \|^2]^{1/2} < \infty.
$$

Then the series $\Gamma = \sum_{j\in\Z} Cov(y_0,y_j)$ converges (coordinatewise) absolutely and the we have
$$
	\Big\{ [0,1]\ni t\mapsto \frac{1}{\sqrt{n}} \sum_{j=1}^{[nt]} y_j \Big\} \Rightarrow W_\Gamma \qquad (n\to \infty),
$$
where the weak convergence is in the $d$-dimensional Skorohod space $D_d[0,1]$ and $W_\Gamma$ denotes a mean zero Gaussian process with covariance $Cov( W_\Gamma(s),W_\Gamma(t)) = \Gamma \cdot \min\{s,t\}$.
\end{theorem}

\begin{theorem}[Lemma 3.2 in Rademacher and Krebs \cite{rademacher2023twosample}]\label{L:USLLN}
Let $(Z_n)_n$ be a sequence of identically distributed, integrable random variables and assume that the corresponding partial sum $\ol{Z}_k = \sum_{i=1}^k Z_i$ satisfies 
a SLLN, i.e., $k^{-1}\ol{Z}_k \to \E[Z_0]$ almost surely.
Then for all $\delta$, $\epsilon>0$ and $\gamma \in (0,1)$ there is a constant $M\in \N$ such that
\begin{align*}
	\p\left( \max_{k\le m } \Big( \frac{k+1}{m}\Big)^\gamma \big| k^{-1} \ol Z_k - \E[Z_0]\big| > \epsilon \right) < \delta \qquad \forall\, m\ge  M.
\end{align*}
\end{theorem}

\begin{theorem}[Corollary 3.3 in Rademacher and Krebs \cite{rademacher2023twosample}]\label{L:SLLN}
Any stationary, $L^2$-$m$-approximable sequence $(Z_n)_n$ satisfies the strong law of large numbers, that is $n^{-1}(Z_1+\ldots+Z_n) \to \E[Z_0]$ almost surely.
\end{theorem}

\section*{Acknowledgments}
Johannes Krebs and Daniel Rademacher are grateful for the financial support of the German Research Foundation (DFG), Grant Number KR-4977/2-1.


\end{document}